\documentclass[11pt,a4paper]{article}
\usepackage{epsf, graphicx}
\usepackage{latexsym,amsfonts,amsbsy,amssymb}
\usepackage{amsmath,amsthm}
\usepackage{mathrsfs}
\usepackage{cite}
\usepackage{cleveref}
\usepackage{enumerate}
\usepackage{enumitem}
\usepackage{color}
\usepackage{xcolor}   
\usepackage{tikz}
\usetikzlibrary{matrix,chains}

\usepackage{tikz}\usetikzlibrary{positioning} 

\usetikzlibrary{calc}\usetikzlibrary{intersections,through}

\makeatletter

\@addtoreset{equation}{section} \makeatother
\makeatletter \@addtoreset{equation}{section} \makeatother
\setlist[enumerate,1]{label=(\arabic*).,font=\textup,
leftmargin=7mm,labelsep=1.5mm,topsep=0mm,itemsep=-0.8mm}
\setlist[enumerate,2]{label=(\alph*).,font=\textup,
leftmargin=7mm,labelsep=1.5mm,topsep=-0.8mm,itemsep=-0.8mm}

\newtheorem{thm}{Theorem}[section]
\newtheorem*{conjC}{Conjecture}

\newtheorem{lem}[thm]{Lemma}
\newtheorem{conj}[thm]{Conjecture}

\newtheorem{defi}[thm]{Definition}
\newtheorem{prop}[thm]{Proposition}
\newtheorem{coro}[thm]{Corollary}
\newtheorem{hyp}[thm]{Hypothesis}

\newcommand{\IBr}{{\rm IBr}}
\newcommand{\IBrd}{{\rm IBrd}}
\newcommand{\which}{\,|\,}
\newcommand{\bwhich}{\,\big|\,}
\newcommand{\Bwhich}{\,\Big|\,}
\newcommand{\Irr}{{\rm Irr}}
\newcommand{\bl}{{\rm bl}}
\newcommand{\Tr}{{\rm Tr}}
\newcommand{\Aut}{{\rm Aut}}
\newcommand{\hH}{\mathcal{H}}

\newcommand{\N}{{\rm\bf N}}
\newcommand{\C}{{\rm\bf C}}
\newcommand{\Z}{{\rm\bf Z}}
\newcommand{\Ind}{{\rm Ind}}
\newcommand{\M}{{\rm M}}
\newcommand{\ra}{\rightarrow}
\newcommand{\mK}{\mathcal{K}}
\newcommand{\mO}{\mathcal{O}}
\newcommand{\GL}{{\rm GL}}
\newcommand{\semi}{\rtimes}
\newcommand{\Pj}{\mathcal{P}}
\newcommand{\PQ}{\mathcal{Q}}
\newcommand{\PE}{\mathcal{E}}
\newcommand{\mS}{\mathcal{S}}
\newcommand{\I}{{\rm I}}
\newcommand{\zg}{\unlhd}

\newcommand{\Fp}{\mathbb{F}_p}

\newcommand{\dz}{{\rm dz}}
\newcommand{\rw}{\mathcal{RW}}
\newcommand{\W}{\mathcal{W}}
\newcommand{\Rad}{{\rm Rad}}

\newcommand{\bS}{{\bf S}}
\newcommand{\bO}{{\bf O}}
\newcommand{\diag}{{\rm diag}}
\newcommand{\Br}{{\rm Br}}

\newcommand{\Cl}{{\mathfrak{Cl}}}

\newcommand{\Bl}{{\rm Bl}}

\newcommand{\cW}{\mathscr{W}}

\begin{document}
\title{\bf A reduction theorem for the blockwise Navarro Alperin weight conjecture via $\hH$-triples}
\date{}
 \author{
  Zhicheng Feng \footnote{Shenzhen International Center for Mathematics, Southern University of Science and Technology, Shenzhen 518055, China. E-mail: fengzc@sustech.edu.cn}\quad ,
  Qulei Fu \footnote{Shenzhen International Center for Mathematics, Southern University of Science and Technology, Shenzhen 518055, China. E-mail: quleifu@outlook.com}\quad ,
 Yuanyang Zhou  \footnote{School of Mathematics and Statistics, Cental China Normal University, Wuhan 430079, China. E-mail: zhouyuanyang@ccnu.edu.cn}
 }
 \maketitle
\begin{abstract}
The Navarro Alperin weight (NAW) conjecture (also known as the Galois Alperin weight conjecture) has been reduced to the inductive NAW condition for simple groups.  
We proceed in two steps to refine this reduction.  
First, we propose the \emph{blockwise Navarro Alperin weight (BNAW) conjecture} and define its associated \emph{inductive BNAW condition}.  
Second, assuming the inductive NAW (respectively, BNAW) condition for simple groups, we establish a stronger version of the NAW (respectively, BNAW) conjecture in terms of central (respectively, block) isomorphism of $\hH$-triples.

\end{abstract}

\section{Introduction}
The \emph{Alperin weight (AW)} conjecture, proposed by Alperin in 1987 \cite{Al87}, has been a central problem in the modular representation theory. Together with the earlier \emph{McKay conjecture} (formulated by McKay in 1971 \cite{Mc71,Mc72}), it belongs to the class of \emph{global–local counting conjectures}, which relate the representation theory of a finite group to that of its local subgroups.
The work of Navarro and Tiep \cite{NT11} marked a breakthrough by \emph{reducing} the AW conjecture to simple groups. They proved that the conjecture holds provided every finite non-abelian simple group satisfies the so-called \emph{inductive AW condition} (referred to as the AWC-good condition in \cite{NT11}), a requirement substantially stronger than the original conjecture.
In 2013, Sp\"ath \cite{Spa13} also achieved a reduction for the \emph{blockwise Alperin weight (BAW)} conjecture.
For progress on the inductive investigation of the Alperin weight conjecture, we refer to the survey paper \cite{FZ22}.

Meanwhile, investigations concerning the action of Galois automorphisms on these conjectures have been progressing. It began with Navarro’s work \cite{Na04} in 2004, which strengthened the McKay conjecture by incorporating Galois automorphisms.
Turull then \cite{Tu14} proved this enhanced AW conjecture in the case of $p$-solvable groups. 
Our earlier work \cite{FFZ} introduced the \emph{Navarro Alperin weight (NAW)} conjecture (also referred to as the \emph{Galois Alperin weight conjecture}), which extends this perspective by considering both Galois automorphisms and group automorphisms. 

Although some progress has been made on studying the action of Galois automorphisms on the AW bijection \cite{Tu14,Duhu25}, no reduction theorem for any Galois version of the AW conjecture was available prior to \cite{FFZ}. In that work, the same authors reduced the NAW conjecture to simple groups and verified the inductive NAW condition for simple Lie groups at their defining characteristic.
We believe that following this approach will reveal a deeper understanding of the AW bijection and ultimately lead to a complete resolution to these conjectures. At the same time, it brings new and significant challenges in the explicit calculation of characters of finite groups of Lie type.

In this paper, we investigate the \emph{blockwise Navarro Alperin weight (BNAW)} conjecture (see Conjecture~\ref{ConjBGAW}) and its reduction to simple groups.
We begin with the necessary notation.

Throughout this paper, we fix a prime $p$ and a $p$-modular system  $(\mK,\mO,F)$ that is sufficiently large for all finite groups under consideration, which means all  group algebras in this paper split over $\mK$ and $F.$
We assume that the residue field $F$ is finite.
Let $\hH \leq \Aut(\mK)$ be a finite abelian group of automorphisms of $\mK$ preserving $\mO$, which therefore induces automorphisms on the residue field $F = \mO/\mathfrak{p},$ where $\mathfrak{p}$ is the unique maximal ideal of $\mO$. 
We assume $\hH$ surjects onto $\Aut(F)$ under the canonical projection $\Aut(\mK, \mO) \to \Aut(F)$ (see \cite[Section 1]{FFZ}).

Let $G$ be a finite group.
A complex irreducible character $\chi\in\Irr(G)$  is said to have \emph{$p$-defect zero} if $\chi(1)_p = |G|_p$, where $n_p$ denotes the $p$-part of an integer $n$. 
Such characters restrict irreducibly to Brauer characters and correspond to projective irreducible modules, hence can be regarded as projective irreducible Brauer characters of $G$ (see the next section).
A \emph{$p$-weight} of $G$ is a pair $(Q,\varphi)$ where $Q$ is a $p$-subgroup of $G,$ and $\varphi \in \Irr(\N_G(Q)/Q)$ has $p$-defect zero.
As developed in Section \ref{sec.Not.&Pre.}, such $\varphi$ corresponds to an  irreducible $F\N_G(Q)$-module with vertex $Q.$
Denote by $\W(G)/\sim_G$ the set of $G$-conjugacy classes of the $p$-weights of $G.$

For a finite group $G,$ the NAW conjecture from \cite{FFZ} asserts the existence of an $\hH \times \Aut(G)$-equivariant bijection between the sets $\IBr(G)$ and $\W(G)/\sim_G$.  
The BNAW conjecture refines this by requiring that the corresponding characters under this bijection lie in  Brauer corresponding blocks. For $\theta\in \IBr(G),$ we denote by $\bl(\theta)$  the block of $G$ containing $\theta,$ and see Section \ref{sec.Not.&Pre.} for the precise definition and notation of block induction. 

\begin{conj}[blockwise Navarro Alperin weight conjecture]\label{ConjBGAW}
    For any finite group $G,$ there exists an $\hH\times\Aut(G)$-equivariant bijection \[\Omega:\IBr(G)\ra \W(G)/\sim_G\] such that $\bl(\varphi)^G=\bl(\theta)$ for any $\theta\in\IBr(G)$ and $(Q,\varphi)\in \Omega(\theta).$
\end{conj}

Where $\varphi$ is regarded as a character of $\N_G(Q).$

In this paper, we introduce the inductive BNAW condition and reduce the  BNAW conjecture to the case of simple groups. 
In fact, we establish a self-reduction of the BNAW (resp. NAW) conjecture by proving that the inductive BNAW (resp. NAW) condition holds for all finite groups if and only if it holds for all quasi-simple groups.
It is worth noting that the analogous problem for the McKay conjecture—namely, the reduction theorem for the Alperin–McKay–Navarro conjecture (Conjecture B in \cite{Na04})—remains open. Our results may provide new insights toward resolving this question.

We now present the strengthened NAW  conjecture and strengthened BNAW conjecture  via the corresponding inductive condition, which are respectively much stronger than the original  conjectures.
The framework of $\hH$-triples with the associated central isomorphisms (resp. block isomorphisms)  currently provides the most efficient characterization of the corresponding inductive conditions. Originally introduced by Navarro, Sp\"ath and Vallejo \cite{NSV20}, this formalism will be properly defined within our current setting in Section \ref{sec.Not.&Pre.}.

\begin{conj}[Inductive NAW conjecture]\label{ConjA}
  Let $G\zg A$ be finite groups. Then there exists an $\hH\times A$-equivariant bijection $\Omega:\IBr(G)\ra\W(G)/\sim_G$ such that $$(A_{\theta^{\hH}},G,\theta)_{\hH}\geqslant_c (\N_A(Q)_{\varphi^{\hH}},\N_G(Q),\varphi)_{\hH}$$ for any $\theta\in\IBr(G)$ and $(Q,\varphi)\in\Omega(\theta).$
\end{conj}

\begin{conj}[Inductive BNAW conjecture]\label{ConjB}
  Let $G\zg A$ be finite groups. Then there exists an $\hH\times A$-equivariant bijection $\Omega:\IBr(G)\ra\W(G)/\sim_G$ such that $$(A_{\theta^{\hH}},G,\theta)_{\hH}\geqslant_b (\N_A(Q)_{\varphi^{\hH}},\N_G(Q),\varphi)_{\hH}$$ for any $\theta\in\IBr(G)$ and $(Q,\varphi)\in\Omega(\theta).$
\end{conj}

Where $\varphi$ is regarded as a Brauer character of $\N_G(Q)$, see the remark in Section 2.

We say that \emph{Conjecture \ref{ConjA} (resp.\ Conjecture \ref{ConjB}) holds for $G$ at the prime $p$} if it holds for every choice of $G \zg A$. Clearly, Conjecture \ref{ConjB} implies Conjecture \ref{ConjA}. Note also that the condition 
\[
(A_{\theta^{\hH}},G,\theta)_{\hH}\geqslant_b (\N_A(Q)_{\varphi^{\hH}},\N_G(Q),\varphi)_{\hH}
\]
in Conjecture \ref{ConjB} implies $\bl(\varphi)^G = \bl(\theta)$.

It is implicit in the proof of Navarro and Tiep \cite{NT11} that, for an arbitrary finite group, the AW bijection is equivariant under the action of group automorphisms, provided the inductive condition holds for all simple groups.
Recent work by Martínez, Rizzo, and Rossi \cite{MRR23} strengthens the original reductions by Navarro--Tiep and Sp\"ath \cite{NT11,Spa13}, showing that the inductive AW (resp. BAW)  condition holds for all finite groups if and only if it holds for all quasi-simple groups.
Analogous reduction theorems have been established by Rossi for the McKay conjecture \cite{Ro23} and for Dade’s character triple conjecture \cite{Ro25}.
As noted in \cite{MRR23}, this enhanced reduction theorem of the (inductive) McKay Conjecture from \cite{Ro23} is used in the verification of the inductive McKay condition for specific families of simple groups.
We adapt these formulations to our setting.

Recall that a group $X$ is called  \emph{involved} in $G$ if there exist subgroups $N \zg H \leq G$ with $H/N \cong X$.
In this paper, we prove the following.

\begin{thm}\label{thmC}
  Let $G$ be a finite group and $p$ a prime number. If Conjecture \ref{ConjA} (resp. Conjecture \ref{ConjB}) holds at $p$ for every universal $p'$-covering group of any finite non-abelian simple group involved in $G,$ then Conjecture \ref{ConjA} (resp. Conjecture \ref{ConjB}) holds for $G.$
\end{thm}

Since Conjecture \ref{ConjB} holds trivially for $p'$-groups, and every non-abelian finite simple group involved in a $p$-solvable group is of $p'$-order, the theorem above yields the following corollary. This also gives an alternative proof of a stronger version of \cite[Corollary 7.3]{Tu14}. 

\begin{coro}
Let $G$ be a finite $p$-solvable group. Then Conjecture \ref{ConjB} holds for $G$ at  $p$.
\end{coro}

In this paper, we use the second cohomology class associated with an $\hH$-triple developed in \cite{Fu24}, and construct a new $\hH$-triple whose normal subgroup is central in the stabilizer and which shares the same Clifford theory as the original one (see Theorem \ref{thm:centrali.}).
This construction generalizes \cite[Theorem 8.28]{Nav:char.} and, in particular, enables us to incorporate central and block isomorphisms of $\hH$-triples into the reduction framework of Navarro and Tiep \cite{NT11}.
In addition, we generalize results such as \cite[Theorem 3.14]{NS14} and  \cite[Theorem 5.8]{MRR23} in the setting of $\hH$-triples (see Theorems  \ref{thm:H-tri(ind)} and \ref{thm:DGN-cor}), which play a crucial role in the final reduction.

We say that the \emph{inductive NAW condition} (respectively, \emph{inductive BNAW condition}) holds for a finite non-abelian simple group $L$ if Conjecture \ref{ConjA} (respectively, Conjecture \ref{ConjB}) holds for the universal $p'$-covering group of $L$.


The following result, which generalizes  \cite[Theorem~C]{FFZ}, provides an evidence for the validity of the inductive BNAW condition for the finite simple  groups.

\begin{thm}\label{thmD}
  The inductive BNAW condition holds for any finite non-abelian simple group of Lie type at its defining characteristic.
\end{thm}

In Section 2, we introduce the basic notation and preliminaries, in particular the definition of $\hH$-triples and the isomorphisms between them. 
Section 3 establishes several results that allow us to construct new isomorphisms of $\hH$-triples from previously known ones.
Sections 4–6 contain three main results that are essential to the final reduction.
In Section 4, we give methods for reducing any $\hH$-triple to one whose normal subgroup is central in the  stabilizer.
Section 5 shows that isomorphisms of $\hH$-triples behave compatibly with the Clifford correspondence.
Section 6 studies bijections of characters lying above the Dade–Glauberman–Nagao correspondence and proves that the corresponding characters in these bijections satisfy the isomorphims of $\hH$-triples.
In Section 7, we move Conjectures \ref{ConjA} and \ref{ConjB} from quasi-simple groups to central extensions of direct products of isomorphic simple groups. Section 8 then provides the final reduction, that is, the proof of Theorem \ref{thmC}.
As evidence supporting the validity of this approach, we prove Theorem \ref{thmD} in Section 9.

\section{Notation and Preliminary results}\label{sec.Not.&Pre.}

In what follows, we introduce some basic notation and conventions used throughout this paper.

Throughout this paper, we fix $p$ to be a prime and every group we consider is finite. For ordinary and Brauer characters, our notation follows \cite{NagTsu}.

Let $(\mathcal{K},\mathcal{O},F)$ be the $p$-modular system and $\hH$ the group of Galois automorphisms as defined in Section~1 (see \cite[Section 1]{FFZ}) for more details).
Both $\mathcal{K}$ and $F$ are sufficiently large for all finite groups under consideration.
The group $\hH$ is a finite abelian group of automorphisms of $\mathcal{K}$ that stabilize $\mathcal{O}$.
Elements in $\hH$ are composed from the left. 
For elements $x \in \mathcal{O} \cup F$ and $\sigma \in \hH$, we denote the action by $x^\sigma$.

Unless otherwise stated, all characters in this paper are Brauer characters.

Two group representations (or projective representations) $X,X':G\ra \GL_m(F)$ are called \emph{similar}, denoted $X' \sim X$, if there exists a $T\in\GL_m(F)$ such that $X'(g)=TX(g)T^{-1}$ for all $g\in G.$ 

For all projective representations $\rho\colon G \to \GL_m(F)$ considered in this paper, we always assume that 
    \( \rho(1) = I_m \).
Consequently, the associated factor set $\alpha\colon G \times G \to F^\times$ is \emph{normalized}, satisfying:
    \[ \alpha(1,g) = \alpha(g,1) = 1 \quad \text{for all } g \in G. \]

Let $G$ be a finite group.
We denote by $G_{p'}$ the set of $p$-regular elements of $G$. 

An ordinary character $\chi \in \Irr(G)$ is said to have \emph{$p$-defect zero}  if $\chi(1)_p = |G|_p$. The set of all irreducible defect zero characters of $G$ is denoted by $\dz(G)$. 

For any $\chi \in \dz(G)$, the restriction $\chi^\circ := \chi|_{G_{p'}}$ yields an irreducible Brauer character of $G$. 
They can be regarded as projective irreducible Brauer characters of $G$, meaning they correspond to projective irreducible $FG$-modules. Equivalently, they have trivial vertices in the sense of \cite[Chapter IV, Section 3]{NagTsu}.


Let $H \leq G$ be a  subgroup and $\lambda \in \IBr(H).$ 
 We denote by $\IBr(G \which \lambda)$ the set of irreducible Brauer characters of $G$ lying over $\lambda$. 
The  same convention follow for ordinary characters. 
Let $\chi\in\IBr(G),$ we denote by $\IBr(H\which \chi)\subseteq\IBr(H)$ the set of irreducible constituents of $\chi_H.$

For a non-negative integer $a$, we define $\IBr(G \which \lambda, p^a)$ to be the subset of $\IBr(G \which \lambda)$ consisting of characters with a vertex of order exactly $p^a$. Similarly, $\IBr(G \which p^a)$ denotes the subset of $\IBr(G)$ consisting of characters with a vertex of order $p^a$.
By vertices of an irreducible Brauer character of $G,$ we mean the minimal $p$-subgroups $Q \leq G$ for which the corresponding $FG$-module is relatively $Q$-projective. 
For  details on vertices, we refer to \cite[Chapter IV, Section 3]{NagTsu}.

For a $p$-weight $(Q,\varphi)$ of $G$, we have $Q$ is a $p$-radical subgroup of $G,$ which means that $Q=\textbf{O}_p(\N_G(Q)),$ where $\textbf{O}_p(\N_G(Q))$ denotes the maximal normal $p$-subgroup of $\N_G(Q)$ (see \cite[Section 2]{NT11}).
Denote by $\Rad(G)$ the set of $p$-radical subgroups of $G.$
By \cite[Lemma 2.1(2)]{FFZ}, there is a natural bijection
\begin{equation}\label{equ:new2.1}
  \dz(\N_G(Q)/Q)\ra \IBr(\N_G(Q)\which |Q|),\chi\mapsto\widetilde{\chi^{\circ}},
\end{equation}
where $\widetilde{\chi^{\circ}}$ is the inflation of $\chi^{\circ}$ to $\N_G(Q)$. 

\noindent\emph{Remark.
Through this bijection, the set $\W(G)$ of $p$-wights of $G$ can be identified with the set 
\begin{equation}\label{equ:new2.2}
  \W^{\circ}(G):=\{(Q,\varphi)\which Q\in\Rad(G),\varphi\in\IBr(\N_G(Q)\which |Q|).\}
\end{equation}
To be more precise with Conjectures \ref{ConjA} and \ref{ConjB}, the bijection involved is in fact   \[\Omega:\IBr(G)\ra \W^{\circ}(G)/\sim_G.\]
}

We let $\Aut(G)$ be the automorphism group of $G,$ where composition is taken from the left.
 For any  $g \in G$ and $\phi \in \Aut(G)$, we write the action as
\(
g^\phi := \phi(g).
\)

Let $H \leq G$ be a subgroup and $\chi$ an ordinary or Brauer character of $H$. 
For any element $a =( \sigma, g) \in \hH \times G,$  we define $H^a := H^g = g^{-1}Hg$, and character $\chi^a$ of $H^a$ given by $\chi^a(x) = \chi(gxg^{-1})^\sigma$ for all $x \in H^a$.
For $a = (\sigma,\phi) \in \hH \times \Aut(H)$, we define the character $\chi^a$ of $H$ by $\chi^a(x) = \chi(x^{\phi^{-1}})^\sigma$ for all $x \in H$.

Let $\mathcal{P}$ be a projective representation of $H$ of degree $m$. For $a = (\sigma,g) \in \mathcal{H} \times \N_G(H)$ (respectively, $a = (\sigma,\phi) \in \mathcal{H} \times \Aut(H)$), we define the function $\mathcal{P}^a \colon H \to \GL_m(F)$ by:
\[
\mathcal{P}^a(x) = \mathcal{P}(gxg^{-1})^\sigma \quad 
\text{(resp.\ } \mathcal{P}^a(x) = \mathcal{P}(x^{\phi^{-1}})^\sigma \text{)}.
\]
Then $\mathcal{P}^a$ is again a projective representation of $H$. 

We denote by $\Bl(G)$ the set of ($p$-)blocks of $G.$
If $\chi\in\IBr(G),$ we write $\bl(\chi)$ for the block of $G$ containing $\chi.$
If $H\leq G$ and $b\in\Bl(H),$ we write $b^G$ for the induced block of $b$ to $G$ (if it is defined), see \cite[Chapter V, Section 3]{NagTsu} for its definition.
We write $\Cl_G(x)$ for the conjugacy class of $G$ containing $x\in G.$

Let $\mathscr{A}$ be a group acting on the right on sets $X_1$ and $X_2$. A map $f \colon X_1 \to X_2$ is called \emph{$\mathscr{A}$-equivariant} if 
\(
f(x^a) = f(x)^a\) for all $ x \in X_1 $ and $a \in \mathscr{A}.$

For any subgroup $N \leq G$, we denote by $\mathcal{S}(G,N)$ the set of all subgroups of $G$ containing $N$. 

Given a character $\theta \in \IBr(N)$, we write $G_\theta$ for its stabilizer in $G$.
Note that we do not require $N \zg G$, but we certainly have $N \zg G_{\theta}$.
For another subgroup $H \leq G$, we define $\N_G(N, H) := \N_G(N) \cap \N_G(H)$.

Let $Q \zg G$ be a normal subgroup of $G$ and write $\bar{G} = G/Q$. 
We establish some conventions on the bar.
For elements $g \in G$, we write $\bar{g} = gQ$ for their images in $\bar{G}$. For subgroups $K \leq G$, we denote $\bar{K} = KQ/Q$. When $K \geq Q$ and $\theta \in \IBr(K)$ with $Q \leq \ker(\theta)$ (in particular when $Q$ is a $p$-group), we let $\bar{\theta}$ be the Brauer character of $\bar{K}$ satisfying $\bar{\theta}(\bar{x}) = \theta(x)$ for all $x \in K$. This notation will be used consistently throughout.


We now introduce character triples and their partial order relations. The theory is well-developed for both ordinary characters \cite{NS14,Spa17} and modular characters \cite{SV16,MRR23}. Standard references include \cite[Chapter 10]{Navarro:McKay} and \cite[Section 3]{MRR23}.

A \emph{modular character triple} $(G,N,\theta)$ consists of finite groups $N \zg G$ and a $G$-invariant Brauer character $\theta \in \IBr(N)$. A \emph{projective representation associated with $\theta$} is a map $\Pj: G \to \GL_{\theta(1)}(F)$ satisfying:
\begin{enumerate}
\item $\Pj|_N$ affords $\theta,$ and
\item $\Pj(g)\Pj(n) = \Pj(gn)$ and $\Pj(n)\Pj(g) = \Pj(ng)$ for all $n \in N$, $g \in G.$
\end{enumerate}
Let $\bar G=G/N.$
The factor set $\alpha:G\times G\ra F^{\times}$ of $\Pj$ is constant on $N\times N$-cosets, and thus induces  a factor set $\bar\alpha:\bar G\times\bar G\ra F^{\times}$ satisfying \[\bar\alpha(\bar g,\bar h)=\alpha(g,h),\quad\forall g,h\in G.\]

If $X$ is a subset of $G,$ then we define \[\Pj(X^+)=\sum_{x\in X}\Pj(x).\]
We let $\Pj(X^+)$ be the zero matrix if $X$ is empty.

For brevity, the term  ``character triple''  will always refer to its modular version.  This convention holds unless explicitly stated otherwise.

\begin{defi}\label{de:char-tri}
  Let $(G,N,\theta)$ and $(H,M,\varphi)$ be character triples with $G=NH$ and $M=N\cap H.$ 
  We identify $G/N\cong H/M.$
  Let $\bar G=G/N$ and thus $\bar{G}=\bar{H}.$
  \begin{enumerate}
    \item If there exist projective representations $\Pj$ of $G$ and $\Pj'$ of $H$ associated with $\theta$ and $\varphi$ with factor sets $\bar\alpha$ and $\bar\alpha',$ respectively, such that $\bar{\alpha}=\bar\alpha',$ then we write $$(G,N,\theta)\geqslant (H,M,\varphi).\qquad \text{(Isomorphism of character triples.)}$$ 
    \item Under condition (1), if additionally $\C_G(N)\subseteq H,$ and for any $c\in\C_G(N),$ the scalar matrices $\Pj(c)$ and $\Pj'(c)$ are associated with the same scalar, then we write $$(G,N,\theta)\geqslant_c (H,M,\varphi).\qquad \text{(Central isomorphism of character triples.)}$$ 
    \item Under conditions (1)-(2), if furthermore  for any $x\in G$ the scalar matrices $\Pj\big(\Cl_{\langle N,x\rangle}(x)^{+}\big)$ and $\Pj'\big((\Cl_{\langle N,x\rangle}(x)\cap H)^{+}\big)$ are associated with the same scalar, 
        then we write $$(G,N,\theta)\geqslant_b (H,M,\varphi).\qquad \text{(Block isomorphism of character triples.)}$$ 
  \end{enumerate}
\end{defi}

For brevity, we write
\begin{equation}\label{equ:Basic1}
  (G,N,\theta)\geqslant_{\ast} (H,M,\varphi),
\end{equation} where $\ast\in\{\emptyset,c,b\}$ corresponds to conditions (1)-(3) in the definition, respectively. In these cases, we say  $(\Pj,\Pj')$ \emph{gives} (\ref{equ:Basic1}).

Suppose that $(\mathcal{P}, \mathcal{P}')$ gives
\((G,N,\theta) \geqslant (H,M,\varphi),\)
 and write $\bar G=G/N.$
By \cite[Theorem 10.13]{Navarro:McKay}, for any $J \in \mathcal{S}(G,N)$, the pair $(\mathcal{P}_J, \mathcal{P}'_{J \cap H})$ induces a bijection
\begin{equation*}
  \nu_J \colon \IBr(J \which \theta) \to \IBr(J \cap H \which \varphi),
\end{equation*}
such that if $\mathcal{Q} \otimes \mathcal{P}_J$ affords $\chi \in \IBr(J \which \theta)$, then $\mathcal{Q} \otimes \mathcal{P}'_{J \cap H}$ affords $\nu_J(\chi)$, where $\mathcal{Q}$ is a projective representation of $\bar{J} = \overline{J \cap H}.$ 
We say that $\nu := \{\nu_J \which J \in \mathcal{S}(G,N)\}$ is the \emph{isomorphism of character triples corresponding to $(\mathcal{P}, \mathcal{P}')$.}

The following theorem is well known (see, for example, \cite[Theorem 10.16]{Navarro:McKay} and \cite[Lemma 3.5]{MRR23}).
We state it here for the reader's convenience.

\begin{thm}\label{thm:2.two.equiv.}
We have the following two equivalences.
\begin{enumerate}
  \item Suppose that $(\mathcal{P}, \mathcal{P}')$ gives $(G, N, \theta) \geqslant (H, M, \varphi)$, and $\nu$ is the isomorphism of character triples corresponding to $(\mathcal{P}, \mathcal{P}')$. 
  Then the following conditions are equivalent:
  \begin{enumerate}
    \item $(\mathcal{P}, \mathcal{P}')$ gives $(G, N, \theta) \geqslant_c (H, M, \varphi)$.
    \item $\IBr(\C_J(N) \which \chi) = \IBr(\C_J(N) \which \nu_J(\chi))$ for all $J \in \mathcal{S}(G, N)$ and $\chi \in \IBr(J \which \theta)$. 
  \end{enumerate}

  \item Suppose that $(\mathcal{P}, \mathcal{P}')$ gives $(G, N, \theta) \geqslant_c (H, M, \varphi)$, and $\nu$ is the  isomorphism of character triples corresponding to $(\mathcal{P}, \mathcal{P}')$. 
  Then the following conditions are equivalent:
  \begin{enumerate}
    \item $(\mathcal{P}, \mathcal{P}')$ gives $(G, N, \theta) \geqslant_b (H, M, \varphi)$.
    \item There exists a defect group $D$ of $\bl(\varphi)$ such that $\C_G(D) \leq H$ and 
    \[
      \bl(\chi) = \bl(\nu_J(\chi))^J
    \]
    for all $J \in \mathcal{S}(G, N)$ and $\chi \in \IBr(J \which \theta)$.
  \end{enumerate}
\end{enumerate}
\end{thm}

In order to consider Galois automorphisms, we work with $\mathcal{H}$-triples as defined in \cite{NSV20}. 
Let $N \leq G$ be finite groups and $\theta \in \IBr(N)$. 
We denote by $\theta^{\mathcal{H}}$ the $\mathcal{H}$-orbit of $\theta$ and by $G_{\theta^{\mathcal{H}}}$ the stabilizer of $\theta^{\mathcal{H}}$ in $G$. 
Explicitly,
\[
G_{\theta^{\mathcal{H}}} = \bigl\{ g \in G \bigm| \text{there exists } \sigma \in \mathcal{H} \text{ such that } \theta^{g} = \theta^{\sigma} \bigr\}.
\]
Note that $N\zg G_{\theta^{\hH}}$  and $G_{\theta} \zg G_{\theta^{\mathcal{H}}}$. 
When $G = G_{\theta^{\mathcal{H}}}$, we write $(G, N, \theta)_{\mathcal{H}}$ and call it a \emph{modular $\mathcal{H}$-triple} (or simply an \emph{$\mathcal{H}$-triple}).

Let $\mathcal{P}$ be a projective representation of $G_{\theta}$ associated with $\theta$. 
For any $a \in (\mathcal{H} \times G)_{\theta}$, there exists a unique function $\mu_a \colon G_{\theta} \to F^{\times}$ such that $\mathcal{P}^a \sim \mu_a \mathcal{P}$, where the projective representation $\mu_a\mathcal{P} \colon G_{\theta} \to \GL_{\theta(1)}(F)$ is defined by
\[
(\mu_a\mathcal{P})(g) = \mu_a(g)\mathcal{P}(g) \quad \text{for all } g \in G_{\theta}.
\]
The function $\mu_a$ is constant on $N$-cosets, and therefore may be viewed as a map $\mu_a \colon G_{\theta}/N \to F^{\times}.$ (See \cite[Remark 1.3]{NSV20}). 

We observe that $\mu_a$ depends only on the $N$-coset of $a$ in $(\mathcal{H} \times G)_{\theta}$.

\begin{lem}\label{lem:H-tri.ind}
Let $(G,N,\theta)_{\mathcal{H}}$ be an $\mathcal{H}$-triple and $\mathcal{P}$ a projective representation of $G_\theta$ associated with $\theta$. For any $a \in (\mathcal{H} \times G)_\theta$, let $\mu_a: G_\theta \to F^\times$ be the unique function satisfying $\mathcal{P}^a \sim \mu_a \mathcal{P}$. Then for all $n \in N$, we have 
\[ \mu_{an} = \mu_a. \]
\end{lem}
\begin{proof}
Fix $n \in N$. For any $x \in G_\theta$, we compute
\begin{align*}
 \Pj^{an}(x)&=(\Pj^a)^n(x)=\Pj^a(nxn^{-1})\sim \mu_a (nxn^{-1}) \Pj(nxn^{-1})  \\
   &=\mu_a(x)\Pj(n)\Pj(x)\Pj(n)^{-1}\sim\mu_a(x)\Pj(x).
\end{align*}
We conclude that $\mu_{an}(x) = \mu_a(x)$ for all $x \in G_\theta$.
\end{proof}

Let $G$ be a subgroup of a finite group $A$, and let $(G,N,\theta)_{\mathcal{H}}$ be an $\mathcal{H}$-triple. Fix a projective representation $\mathcal{P}$ of $G_\theta$ associated with $\theta$. 
Suppose $K \zg A$ satisfies $K \cap G \leq N$, and let $\bar{A} = A/K$. We establish the following conventions for the bar notation.
First, the factor set $\bar{\alpha}$ of $\mathcal{P}$ is defined as the map $\bar{\alpha} \colon \overline{G_\theta} \times \overline{G_\theta} \to F^\times$ satisfying
\[
\mathcal{P}(g)\mathcal{P}(h) = \bar{\alpha}(\bar{g},\bar{h})\mathcal{P}(gh) \quad \text{for all } g,h \in G_\theta.
\]
Second, we define $(\mathcal{H} \times \bar{G})_\theta$ to be the image of $(\mathcal{H} \times G)_\theta$ under the natural surjection $\mathcal{H} \times G \to \mathcal{H} \times \bar{G}$, thus $(\mathcal{H} \times \bar{G})_\theta\cong (\mathcal{H} \times G)_\theta/(N\cap K).$
Third, for any $a \in (\mathcal{H} \times \bar{G})_\theta$, the map $\mu_a \colon \overline{G_\theta} \to F^\times$ satisfying $\mathcal{P}^a \sim \mu_a \mathcal{P}$ is well-defined by the previous lemma. 
For any  $\chi \in \IBr(J \which \theta)$ where $J \in \mathcal{S}(G,N)$ and any $a \in \mathcal{H} \times \bar{G}$, the  character $\chi^a\in\IBr(J^a\which \theta)$ is well-defined. 

We now introduce the corresponding isomorphisms of $\mathcal{H}$-triples.

\begin{defi}
  Suppose that $(G,N,\theta)_{\hH}$ and $(H,M,\varphi)_{\hH}$ are $\hH$-triples with $G=NH$ and $M=N\cap H.$
  Let $\bar{G}=G/N.$
  Assume that $(\hH\times \bar H)_{\theta}=(\hH\times\bar H)_{\varphi}.$
  Let $\Pj$ and $\Pj'$ be projective representations of $G_{\theta}$ and $H_{\varphi}$ associated with $\theta$ and $\varphi,$ respectively. 
  We assume that $\mu_a=\mu'_a$ for any $a\in(\hH\times\bar H)_{\theta},$ where $\mu_a,\mu'_a:\overline{H_{\varphi}}\ra F^{\times}$ are determined by $\Pj^a\sim\mu_a\Pj,\Pj'^a\sim\mu'_a\Pj',$ respectively.
  Then we write $$(G,N,\theta)_{\hH}\geqslant_{\ast}(H,M,\varphi)_{\hH}$$
  if $(\Pj,\Pj')$ gives $(G_{\theta},N,\theta)\geqslant_{\ast} (H_{\varphi},M,\varphi),$ where $\ast\in\{\emptyset,c,b\}.$
\end{defi}

In the  above setting, we say that $(\Pj,\Pj')$ \emph{gives} 
\begin{equation}\label{equ:prel2}
  (G,N,\theta)_{\hH}\geqslant_{\ast} (H,M,\varphi)_{\hH}.
\end{equation}
Let $\nu$ be the isomorphism of character triples corresponding to $(\Pj,\Pj').$ 
By the arguement in \cite[Lemma 5.6]{FFZ},  for any $J\in\mS(G_{\theta},N)$ the bijection $$\nu_J:\IBr(J\which\theta)\ra\IBr(J\cap H\which\varphi)$$ is $\big(\hH\times\N_H(J)\big)_{\theta}$-equivariant.

The following result, known as the \emph{Butterfly Theorem}, reveals that the isomorphisms of $\mathcal{H}$-triples only depends on the automorphisms of the normal subgroup induced via conjugation by the overgroup.  This theorem was first formally established in \cite[Theorem 5.3]{Spa17}, and we enhance it with Galois automorphisms and block induction.

\begin{thm}\label{thm:butterfly}
 Let $(G,N,\theta)_{\hH}$ and $(H,M,\varphi)_{\hH}$ be $\hH$-triples such that $(G,N,\theta)_{\hH}\geqslant_{\ast}(H,M,\varphi)_{\hH},$ where $\ast\in\{c,b\}.$ Let $(\widehat{G},N,\theta)_{\hH}$ and $(\widehat{H},M,\varphi)_{\hH}$ be $\hH$-triples with $\widehat{G}=N\widehat{H}, M=N\cap\widehat{H},\C_{\widehat{G}}(N)\subseteq \widehat{H}$ and $(\hH\times\widehat{H})_{\theta}=(\hH\times \widehat{H})_{\varphi}.$ 
  Let $\epsilon \colon H \to \Aut(N)$ and $\hat{\epsilon} \colon \widehat{H} \to \Aut(N)$ be the homomorphisms induced by conjugation. If $\epsilon(H) = \hat{\epsilon}(\widehat{H})$, then
\[
(\widehat{G},N,\theta)_{\mathcal{H}} \geqslant_\ast (\widehat{H},M,\varphi)_{\mathcal{H}}.
\]

\end{thm} 
\begin{proof}
The case when $\ast = c$ was established in \cite[Theorem 5.2]{Fu24}. For $\ast = b$, we adapt the proof strategy from \cite[Theorem 4.3]{Spa17}, providing full details here for completeness.

Assume $\ast = b$ and maintain all notation from the proof of \cite[Theorem 5.2]{Fu24}, replacing $\eta$ and $\eta'$ with $\theta$ and $\varphi$ respectively. We have already shown that $(\widehat{\mathcal{P}}, \widehat{\mathcal{P}}')$ gives
\[
(\widehat{G}, N, \theta)_{\mathcal{H}} \geqslant_c (\widehat{H}, M, \varphi)_{\mathcal{H}}.
\]
To complete the proof, it remains to verify that for any $\hat{x} \in \widehat{G}_\theta$ and $\widehat{J} = \langle N, \hat{x} \rangle$, the scalar matrices
\[
\widehat{\mathcal{P}}\big(\Cl_{\widehat{J}}(\hat{x})^{+}\big) \quad \text{and} \quad \widehat{\mathcal{P}}'\big((\Cl_{\widehat{J}}(\hat{x}) \cap \widehat{H})^{+}\big)
\]
are associated with the same scalar.

Choose $\hat x\in\widehat{G}_{\theta}$ and let $\widehat{J}=\langle N,\hat x\rangle.$
 Let $\hat x=\hat t n\hat c,$ where $t\in\mathcal{T},n\in N,\hat c\in\C_{\widehat{G}}(N),$ and let $x=tn\in G_{\theta}.$ 
 Observe that $\hat{x}$ and $x$ induce identical automorphisms on $N$ via conjugation, since $\hat{c} \in \C_{\widehat{G}}(N)$ acts trivially.
 Define the map $\mathcal{L}_{x}:N\ra N,n\mapsto n^{-1}n^{ x^{-1}}.$ 
 Note that $\Cl_{\widehat{J}}(\hat x)=\mathcal{L}_{{x}}(N)\hat x$ and $\Cl_{J}( x)=\mathcal{L}_{x}(N) x,$ where $J=\langle N,x\rangle.$
We compute that
\begin{align*}
\widehat{\mathcal{P}}\big(\Cl_{\widehat{J}}(\hat{x})^{+}\big) 
    &= \sum_{l\in\mathcal{L}_{x}(N)} \widehat{\mathcal{P}}(l\hat{x}) 
    = \sum_{l\in\mathcal{L}_{x}(N)} \mathcal{P}(l)\widehat{\mathcal{P}}(\hat{x}) \\
    &= \sum_{l\in\mathcal{L}_{x}(N)} \mathcal{P}(l)\mathcal{P}(tn)\hat{\mu}(\hat{c})
    = \sum_{l\in\mathcal{L}_{x}(N)} \mathcal{P}(lx)\hat{\mu}(\hat{c}) \\
    &= \mathcal{P}\big(\Cl_{J}(x)^{+}\big)\hat{\mu}(\hat{c}).
\end{align*}
The parallel computation for $\widehat{\mathcal{P}}'$ yields 
 \begin{align*}
   \widehat{\Pj}'\big((\Cl_{\widehat{J}}(\hat x)\cap \widehat{H})^{+}\big)&= \sum_{\substack {l\in\mathcal{L}_{ x}(N)\\l\hat x\in\widehat{H}}}\widehat{\Pj}'(l\hat x)=\sum_{\substack {l\in\mathcal{L}_{ x}(N)\\l\hat x\in\widehat{H}}}\widehat{\Pj}'(\hat xl^{\hat x}) \\
    &= \sum_{\substack {l\in\mathcal{L}_{ x}(N)\\l\hat x\in\widehat{H}}}\Pj'(tnl^{ x})\hat{\mu}(\hat c)
    =\sum_{\substack {l\in\mathcal{L}_{ x}(N)\\l x\in H}}\Pj'(lx)\hat{\mu}(\hat c)\\
    &=\Pj'\big( (\Cl_J(x)\cap H)^{+} \big)\hat\mu(\hat c).
 \end{align*}
Note that for any $l \in \mathcal{L}_x(N)$, we have that $l\hat{x} \in \widehat{H}$ if and only if $lx \in H$, as both conditions are determined by $\hat\epsilon(l\hat{x}) = \epsilon(lx) \in \epsilon(H)$. 
Since by assumption $\mathcal{P}\big(\Cl_{J}(x)^{+}\big)$ and $\mathcal{P}'\big((\Cl_J(x) \cap H)^{+}\big)$ are associated with the same scalar, this completes the proof of the theorem.
\end{proof}

Let $(G,N,\theta)_{\mathcal{H}}$ and $(H,M,\varphi)_{\mathcal{H}}$ be $\mathcal{H}$-triples satisfying 
\begin{equation}\label{equ:pre2}
  (G,N,\theta)_{\mathcal{H}} \geqslant_\ast (H,M,\varphi)_{\mathcal{H}}
\end{equation} for $\ast \in \{c,b\}$. By definition, for any $J \in \mathcal{S}(G,N)$, (\ref{equ:pre2}) restricts to
\[
(J,N,\theta)_{\mathcal{H}} \geqslant_\ast (J \cap H,N,\varphi)_{\mathcal{H}}.
\]
Combining this observation with Theorem \ref{thm:butterfly}, we see that when $H$ is sufficiently large to realize all possible automorphisms of $N$, each pair of $\mathcal{H}$-triples lying above $\theta$ and $\varphi$ will satisfy the isomorphims relation.

Now return to Conjectures \ref{ConjA} and \ref{ConjB}.
They first require an equivariant bijection, say the $\Omega$ there.
We will see the that some assumptions in the above theorem can be easily deduced by this equivariance properties. 
In fact, these are group-theoretical consequences.

\begin{prop}\label{prop:lieover2}
Let $G \zg A$ be finite groups with an $\hH \times A$-equivariant bijection 
$\Omega \colon \IBr(G) \to \W^{\circ}(G)/\sim_G$. For $\theta \in \IBr(G)$ and $(Q,\varphi) \in \Omega(\theta)$, we have:
\begin{enumerate}
    \item $\big(\hH \times \N_A(Q)\big)_\theta = \big(\hH \times \N_A(Q)\big)_\varphi,$ from which we deduce both $\N_A(Q)_\theta = \N_A(Q)_\varphi$ and $\N_A(Q)_{\theta^{\mathcal{H}}} = \N_A(Q)_{\varphi^{\mathcal{H}}}$.
    \item $A_{\theta^{\hH}} = G \N_A(Q)_{\varphi^{\hH}}$ and $G \cap \N_A(Q)_{\varphi^{\hH}} = \N_G(Q)$. 
\end{enumerate}  
\end{prop}

\begin{proof}
Recall that the action of $\mathcal{H} \times A$ on $\mathcal{W}^{\circ}(G)/\sim_G$ is induced by
\[
(Q,\delta)^a = (Q^a, \delta^a) \quad \text{for } (Q,\delta) \in \mathcal{W}(G) \text{ and } a \in \mathcal{H} \times A.
\]
The first equality follows directly from the equivariance of $\Omega$.
For the second statement, given $x \in A_{\theta^{\hH}}$,  there exists $\sigma \in \mathcal{H}$ such that $a = (\sigma, x)$ stabilizes $\theta$. Since $(Q,\varphi)^a$ is $G$-conjugate to $(Q,\varphi)$, we find $g \in G$ satisfying $(Q,\varphi)^{ag} = (Q,\varphi)$. This yields $ag =(\sigma, xg)\in \big(\mathcal{H} \times \N_A(Q)\big)_\varphi, $ and thus \[xg \in \N_A(Q)_{\varphi^{\mathcal{H}}}.\]
Consequently, we obtain 
\(
A_{\theta^{\mathcal{H}}} = G \N_A(Q)_{\varphi^{\mathcal{H}}},
\)
with $G \cap \N_A(Q)_{\varphi^{\mathcal{H}}} = \N_G(Q)$. 
\end{proof}

\noindent\emph{Remark.  In the proposition, the equality 
\[\big(\hH \times \N_A(Q)\big)_\theta = \big(\hH \times \N_A(Q)\big)_\varphi,\]  implies 
\[\big(\hH \times \N_A(Q)_{\theta^{\hH}}\big)_\theta = \big(\hH \times \N_A(Q)_{\theta^{\hH}}\big)_\varphi.\]
Moreover, we have $\C_{A_{\theta^{\hH}}}(N)\subseteq \N_A(Q)_{\varphi^{\hH}}=\N_A(Q)_{\theta^{\hH}} =\N_{A_{\theta^{\hH}}}(Q).$ 
 }

When verifying Conjecture \ref{ConjB} in specific cases, it is often sufficient to work at the block level. 
Let $G$ be a finite group and $B$ a $p$-block of $G$. 
Define
\[
\W^{\circ}(B) := \{(Q,\varphi) \in \W^{\circ}(G) \mid \bl(\varphi)^G = B\},
\]
which is clearly closed under $G$-conjugation.

We say that \emph{Conjecture \ref{ConjB} holds for the block $B$} if for every extension $G \zg A$, the following conditions hold.

\begin{conjC}[Conjecture \ref{ConjB} for a block $B$ of $G$]
  Let $G\zg A$ be finite groups, and consider a block $B$ of $G.$
  Denote by $(\hH\times A)_B$ the stabilizer of $B$ in $\hH\times A.$  Then there exists an $(\hH\times A)_B$-equivariant bijection \[\Omega:\IBr(B)\ra\W^{\circ}(B)/\sim_G\] such that \[(A_{\theta^{\hH}},G,\theta)_{\hH}\geqslant_b (\N_A(Q)_{\varphi^{\hH}},\N_G(Q),\varphi)_{\hH}\] for any $\theta\in\IBr(B)$ and $(Q,\varphi)\in\Omega(\theta).$
\end{conjC}

If Conjecture \ref{ConjB} holds for every $p$-block of $G$, we prove in Section 3 that it then holds for the group $G$ itself.

\section{$\hH$-triples}

In this section, we introduce some properties of the isomorphisms on $\mathcal{H}$-triples, including methods for constructing new isomorphisms between them.

Let $(G,N,\theta)_{\mathcal{H}}$ and $(H,M,\varphi)_{\mathcal{H}}$ be $\mathcal{H}$-triples satisfying
\begin{equation}\label{eq:H-triple-order}
    (G,N,\theta)_{\mathcal{H}} \geqslant_{\ast} (H,M,\varphi)_{\mathcal{H}}
\end{equation}
for $\ast \in \{c,b\}$. Writing $\bar{G} = G/N$, we have $\bar{G} = \bar{H}$.
If the relation \eqref{eq:H-triple-order} is given by $(\mathcal{P}, \mathcal{P}')$, then for any function $\mu \colon \overline{G_{\theta}} \to F^{\times}$ with $\mu(1) = 1$, the same order relation is also given by $(\mu\mathcal{P}, \mu\mathcal{P}')$. This observation yields the following lemma.

\begin{lem}\label{lem:H-tri.tran.}
Let $(G,N,\theta)_{\mathcal{H}}$, $(H,M,\varphi)_{\mathcal{H}}$, and $(L,K,\gamma)_{\mathcal{H}}$ be $\mathcal{H}$-triples satisfying
\[
(G,N,\theta)_{\mathcal{H}} \geqslant_\ast (H,M,\varphi)_{\mathcal{H}}
\quad \text{and} \quad
(H,M,\varphi)_{\mathcal{H}} \geqslant_\ast (L,K,\gamma)_{\mathcal{H}},
\]
where $\ast \in \{c,b\}$. Then the following transitivity holds:
\[
(G,N,\theta)_{\mathcal{H}} \geqslant_\ast (L,K,\gamma)_{\mathcal{H}}.
\]
\end{lem}
\begin{proof}
Suppose that $(G,N,\theta)_{\mathcal{H}} \geqslant_\ast (H,M,\varphi)_{\mathcal{H}}$ is given by $(\mathcal{P},\mathcal{Q})$, and $(H,M,\varphi)_{\mathcal{H}} \geqslant_\ast (L,K,\gamma)_{\mathcal{H}}$ is given by $(\mathcal{Q}',\mathcal{E})$. 

By making suitable adjustments to $(\mathcal{Q}',\mathcal{E})$ (as described in the preceding paragraph), we may assume $\mathcal{Q} = \mathcal{Q}'$. A direct verification shows that $(\mathcal{P},\mathcal{E})$ gives
\[
(G,N,\theta)_{\mathcal{H}} \geqslant_\ast (L,K,\gamma)_{\mathcal{H}}.
\]
The block isomorphism condition from Definition \ref{de:char-tri} is satisfied since the scalar matrices $\mathcal{P}\big(\Cl_{\langle N,x\rangle}(x)^{+}\big)$ and $\mathcal{E}\big((\Cl_{\langle N,x\rangle}(x)\cap L)^{+}\big)$ are associated the scalars for each $x \in G$.
\end{proof}

The following lemma is standard, but we include it here for the reader's convenience.

\begin{lem}\label{lem:XXX}
Let $(G,N,\theta)$ be a character triple and $\mathcal{P}$ a projective representation of $G$ associated with $\theta$.
Let $\mathcal{Q}$ and $\mathcal{E}$ be projective representations of $\overline{G} = G/N$.
  Then \[
\mathcal{Q} \sim \mathcal{E} \quad \text{if and only if} \quad \mathcal{Q} \otimes \mathcal{P} \sim \mathcal{E} \otimes \mathcal{P}.
\]
\end{lem}
\begin{proof}
The ``only if'' part is immediate. For the ``if'' direction, suppose there exists an invertible matrix \( T \) such that for all \( g \in G \),
\[
\mathcal{E}(\bar{g}) \otimes \mathcal{P}(g) = T \left( \mathcal{Q}(\bar{g}) \otimes \mathcal{P}(g) \right) T^{-1}.
\]
Since \( T \) commutes with every \( I_s \otimes \mathcal{P}(n) \), where \( n \in N \) and \( s \) is the degree of \( \mathcal{Q} \) and \( \mathcal{E} \), it follows from \cite[Chapter II, Lemma 4.1(i)]{NagTsu} that \( T = A \otimes I_{\theta(1)} \) for some invertible matrix \( A \). Consequently, \( \mathcal{E} \sim \mathcal{Q} \), completing the proof of the lemma.
\end{proof}
\noindent\emph{Remark. $\mathcal{Q} \otimes \mathcal{P}$ and $\mathcal{E} \otimes \mathcal{P}$ are projective (but not necessarily linear) representations.}

The following lemma will be used in proving Lemma \ref{lem:H-tri.B}.

\begin{lem}\label{lem:H-triA}
Let $(G,N,\theta)_{\mathcal{H}}$ be an $\mathcal{H}$-triple. Set $\bar{G} = G/N$, and let $\mathcal{P}$ be a projective representation of $G_{\theta}$ associated with $\theta$ with factor set $\bar{\alpha}$. For $J \in \mathcal{S}(G_{\theta},N)$ and $\chi \in \IBr(J \which \theta)$, we have:

\begin{enumerate}
    \item Let $\mathcal{Q}$ be a projective representation of $\overline{G_{\chi}}$ with factor set $\bar{\delta}$ satisfying $\mathcal{Q}(1) = I$. Then the following two conditions are equivalent:
    \begin{enumerate}
        \item $\mathcal{Q} \otimes \mathcal{P}_{G_{\chi}}$ is a projective representation of $G_{\chi}$ associated with $\chi$;
        \item The following conditions hold:
            \begin{itemize}
                \item $\bar{\delta}_{\bar{J} \times \bar{J}} = \bar{\alpha}^{-1}_{\bar{J} \times \bar{J}}$;
                \item The representation $\mathcal{Q}_{\bar{J}} \otimes \mathcal{P}_J$ affords $\chi$;
                \item The factor set $\bar{\beta} := \bar{\delta}\bar{\alpha}_{\overline{G_{\chi}} \times \overline{G_{\chi}}}$ is constant on $\bar{J} \times \bar{J}$-cosets.
            \end{itemize}
    \end{enumerate}
    Moreover, such a projective representation $\mathcal{Q}$ exists.

    \item For $a \in (\mathcal{H} \times \bar{G})_{\chi}$ (note that $(\mathcal{H} \times \bar{G})_{\chi} \subseteq (\mathcal{H} \times \bar{G})_{\theta}$), let $\mu_a: \overline{G_{\theta}} \to F^{\times}$ be determined by $\mathcal{P}^a \sim \mu_a \mathcal{P}$. If $\mathcal{Q}$ satisfies the conditions in (1), and $\varrho_a: \overline{G_{\chi}} \to F^{\times}$ is determined by
    \[
    (\mathcal{Q} \otimes \mathcal{P}_{G_{\chi}})^a \sim \varrho_a (\mathcal{Q} \otimes \mathcal{P}_{G_{\chi}}),
    \]
    (where $\varrho_a$ is constant on $\bar{J}$-cosets), 
    then $\varrho_a$ is the unique function $\varrho_a \colon \overline{G_{\chi}} \to F^{\times}$ satisfying
\[
\mathcal{Q}^a \mu_{a,\overline{G_{\chi}}} \sim \varrho_a \mathcal{Q}.
\]
\end{enumerate}
\end{lem}
\begin{proof}
For part (1), the implication (a) $\Rightarrow$ (b) follows immediately from the definitions. Now assume (b) holds.
Since $\bar\beta$ is factor set of $\PQ\otimes\Pj_{G_{\chi}},$
it suffices to verify that
\[
\bar{\beta}(\bar{g},\bar{x}) = \bar{\beta}(\bar{x},\bar{g}) = 1 \quad \text{for all } g \in G_{\chi} \text{ and } x \in J.
\]
By our assumption, we have $\bar{\beta}(\bar{g},\bar{x}) = \bar{\beta}(\bar{g},1) = 1$, and similarly $\bar{\beta}(\bar{x},\bar{g}) = 1$.

To establish the existence of such a $\mathcal{Q}$, let $\mathcal{A}$ be a projective representation of $G_{\chi}$ associated with $\chi$. We may assume $\mathcal{A}_J = \mathcal{E} \otimes \mathcal{P}_J$, where $\mathcal{E}$ is a projective representation of $\bar{J}$ with factor set $\bar{\alpha}^{-1}_{\bar{J} \times \bar{J}}$ such that $\mathcal{E} \otimes \mathcal{P}_J$ affords $\chi$.
Let $m$ be the degree of $\mathcal{E}$. For each $g \in G_{\chi}$, the matrix $\mathcal{A}(g)(I_m \otimes \mathcal{P}(g))^{-1}$ commutes with all $I_m \otimes \mathcal{P}(n)$ for $n \in N$. 
To verify this, observe that for any $n \in N$:
\begin{align*}
\mathcal{A}(n)^{\mathcal{A}(g)} &= \mathcal{A}(n^g), \\
(I_m \otimes \mathcal{P}(n))^{I_m \otimes \mathcal{P}(g)} &= I_m \otimes \mathcal{P}(n^g).
\end{align*}
Since $\mathcal{A}(n) = I_m \otimes \mathcal{P}(n)$ and $\mathcal{A}(n^g) = I_m \otimes \mathcal{P}(n^g)$ for all $n \in N$ and $g \in G_{\chi}$, the commutation relation holds. Applying \cite[Chapter II, Lemma 4.1(i)]{NagTsu}, we obtain a decomposition
\[
\mathcal{A}(g)(I_m \otimes \mathcal{P}(g))^{-1} = \mathcal{Q}(g) \otimes I_{\theta(1)}
\]
for some invertible matrix $\mathcal{Q}(g) \in \mathrm{GL}_m(F)$, and consequently
\[
\mathcal{A}(g) = \mathcal{Q}(g) \otimes \mathcal{P}(g).
\]
The map $\mathcal{Q} \colon G_{\chi} \to \mathrm{GL}_m(F)$ is a projective representation because both $\mathcal{A}$ and $\mathcal{P}$ are projective representations.
Moreover, the relation $\mathcal{A}(g)\mathcal{A}(n) = \mathcal{A}(gn)$ implies $\mathcal{Q}(gn) = \mathcal{Q}(g)$ for all $g\in G_{\chi}$ and $n \in N$, showing that $\mathcal{Q}$ descends to a well-defined projective representation of $\overline{G_{\chi}}$.

For part (2), let $a = (\sigma, \bar{g}) \in \mathcal{H} \times \overline{G}$. Since $\bar{g}$ normalizes $\overline{G_{\chi}}$, the projective representation
\[
\mathcal{Q}^a \mu_{a,\overline{G_{\chi}}} \colon \overline{G_{\chi}} \to \GL_m(F), \quad \bar{x} \mapsto \mathcal{Q}(\bar{g}\bar{x}\bar{g}^{-1})^\sigma \mu_a(\bar{x})
\] 
is well-defined.
By Lemma~\ref{lem:XXX}, for any function $\varrho \colon \overline{G_{\chi}} \to F^\times$, we have:
\[
\mathcal{Q}^a \mu_{a,\overline{G_{\chi}}} \sim \varrho \mathcal{Q} 
\quad \text{if and only if} \quad 
\mathcal{Q}^a \mu_{a,\overline{G_{\chi}}} \otimes \mathcal{P}_{G_{\chi}} \sim \varrho \mathcal{Q} \otimes \mathcal{P}_{G_{\chi}}.
\]
The key observation is that 
\[
(\mathcal{Q} \otimes \mathcal{P}_{G_{\chi}})^a \sim \mathcal{Q}^a \mu_{a,\overline{G_{\chi}}} \otimes \mathcal{P}_{G_{\chi}},
\]
and therefore the uniquely determined function $\varrho_a \colon \overline{G_{\chi}} \to F^\times$ satisfying
\[
(\mathcal{Q} \otimes \mathcal{P}_{G_{\chi}})^a \sim \varrho_a (\mathcal{Q} \otimes \mathcal{P}_{G_{\chi}})
\]
is exactly the one such that
\[
\mathcal{Q}^a \mu_{a,\overline{G_{\chi}}} \otimes \mathcal{P}_{G_{\chi}} \sim \varrho_a \mathcal{Q} \otimes \mathcal{P}_{G_{\chi}}.
\]
This completes the proof of (2).
\end{proof}

We now present several methods for constructing new order relations of $\hH$-triples from existing ones.

\begin{lem}\label{lem:H-tri.B}
Let $(G,N,\theta)_{\mathcal{H}}$ and $(H,M,\varphi)_{\mathcal{H}}$ be $\mathcal{H}$-triples such that $(\mathcal{P},\mathcal{P}')$ gives
\[
(G,N,\theta)_{\mathcal{H}} \geqslant_\ast (H,M,\varphi)_{\mathcal{H}}
\]
for $\ast \in \{\emptyset, c, b\}$. Let $\nu$ be the isomorphism of character triples corresponding to $(\mathcal{P},\mathcal{P}')$. 
For any $J \in \mathcal{S}(G_\theta, N)$ and $\chi \in \IBr(J \which \theta)$, we have $H_{\nu_J(\chi)^{\mathcal{H}}}=H_{\chi^{\hH}}$ and
\begin{equation}\label{equ:H-tri1}
(G_{\chi^{\mathcal{H}}}, J, \chi)_{\mathcal{H}} \geqslant_\ast (H_{\nu_J(\chi)^{\mathcal{H}}}, J \cap H, \nu_J(\chi))_{\mathcal{H}}.
\end{equation}
\end{lem}
\begin{proof}
Set  $\zeta = \nu_J(\chi)$. 
It is straightforward to see that $G_{\chi^{\mathcal{H}}} = JH_{\chi^{\mathcal{H}}}$ and $J\cap H_{\zeta^{\mathcal{H}}}=J\cap H.$
The $\big(\mathcal{H} \times \N_H(J)\big)_{\theta}$-equivariance of $\nu_J: \IBr(J \which \theta) \to \IBr(J\cap H \which \varphi)$ implies $(\mathcal{H} \times H)_{\chi} = (\mathcal{H} \times H)_{\zeta}$, and consequently $H_{\chi^{\mathcal{H}}} = H_{\zeta^{\mathcal{H}}}.$ 
Note that $(\hH\times H)_{\chi}=(\hH\times H_{\chi^{\hH}})_{\chi}$ and  $(\hH\times H)_{\zeta}=(\hH\times H_{\zeta^{\hH}})_{\zeta}.$

Let $\bar G=G/N.$
  Let $\PQ$ be a projective representation of $ \overline{G_{\chi}}$ such that $\PQ\otimes\Pj_{G_{\chi}}$ is a projective representation of $G_{\chi}$ associated with $\chi$ (such $\mathcal{Q}$ exists by Lemma~\ref{lem:H-triA}(1)).
  Since $\overline{G_{\chi}} = \overline{H_{\zeta}}$, Lemma~\ref{lem:H-triA}(1) also shows that $\mathcal{Q} \otimes \mathcal{P}'_{H_{\zeta}}$ is a projective representation of $H_{\zeta}$ associated with $\zeta$, sharing the same factor set as $\mathcal{Q} \otimes \mathcal{P}_{G_{\chi}}$.
  For any $a\in(\hH\times \bar H)_{\chi},$ let the functions $\varrho_a,\varrho'_a:\overline{G_{\chi}}\ra F^{\times}$ be determined by $(\PQ\otimes\Pj_{G_{\chi}})^a\sim \varrho_a(\PQ\otimes\Pj_{G_{\chi}})$ and $(\PQ\otimes \Pj'_{H_{\zeta}})^a=\varrho'_a(\PQ\otimes \Pj'_{H_{\zeta}}),$ respectively.
 Lemma \ref{lem:H-triA}(2) establishes that $\varrho_a=\varrho'_a.$  
This proves that $(\mathcal{Q} \otimes \mathcal{P}_{G_{\chi}}, \mathcal{Q} \otimes \mathcal{P}'_{H_{\zeta}})$ gives \eqref{equ:H-tri1} for $\ast = \emptyset$.

 For the case $\ast = c$, take $c \in \C_{G_{\chi^{\mathcal{H}}}}(J) \subseteq \C_G(N)$. The matrix $\mathcal{Q}(\overline{c}) \otimes \mathcal{P}(c)$ is scalar because $\mathcal{Q} \otimes \mathcal{P}_{G_{\chi}}$ is a projective representation associated with $\chi$.
Therefore $\mathcal{Q}(\overline{c})$ is scalar, so $\mathcal{Q}(\overline{c}) \otimes \mathcal{P}(c)$ and $\mathcal{Q}(\overline{c}) \otimes \mathcal{P}'(c)$ represent the same scalar value.

  Now assume $\ast = b$. For any $x \in G_{\chi} \subseteq G_{\theta}$, let $\Cl_{\langle J,x\rangle}(x) = \bigsqcup_{i=1}^l C_i$ be the decomposition into disjoint $N$-conjugacy classes. By our assumption, there exist scalars $\lambda_i \in F^\times$ such that
\[
\mathcal{P}(C_i^+) = \lambda_i I_{\theta(1)} \quad \text{and} \quad \mathcal{P}'((C_i \cap H)^+) = \lambda_i I_{\varphi(1)}
\]
for each $1 \leq i \leq l$. Choosing representatives $x_i \in C_i$, we compute:
\begin{align}
(\mathcal{Q} \otimes \mathcal{P}_{G_{\chi}})(\Cl_{\langle J,x\rangle}(x)^+) 
= \sum_{i=1}^l \mathcal{Q}(\overline{x}_i) \otimes \mathcal{P}(C_i^+) 
= \left( \sum_{i=1}^l \lambda_i \mathcal{Q}(\overline{x}_i) \right) \otimes I_{\theta(1)}, \label{eq:scalar-sum-G}
\end{align}
and similarly
\begin{align}
(\mathcal{Q} \otimes \mathcal{P}'_{H_{\zeta}})((\Cl_{\langle J,x\rangle}(x) \cap H)^+)
&= \left( \sum_{i=1}^l \lambda_i \mathcal{Q}(\overline{x}_i) \right) \otimes I_{\varphi(1)}. \label{eq:scalar-sum-H}
\end{align}
Since both \eqref{eq:scalar-sum-G} and \eqref{eq:scalar-sum-H} are scalar matrices, they must represent the same scalar value. This completes the proof for this lemma.
\end{proof}

\begin{lem}[Lifting the order relations]\label{lem:H-tri.lift}
Let $(G,N,\theta)_{\mathcal{H}}$ and $(H,M,\varphi)_{\mathcal{H}}$ be $\mathcal{H}$-triples such that $G = NH$ and $M = N \cap H$. Let $Z \subseteq \ker(\theta) \cap \ker(\varphi)$ be a normal subgroup of $G$.
Set $\bar{G} = G/Z$. If $(\bar{G}, \bar{N}, \bar{\theta})_{\mathcal{H}} \geqslant_\ast (\bar{H}, \bar{M}, \bar{\varphi})_{\mathcal{H}}$ for $\ast \in \{\emptyset, c, b\}$, then
\[
(G,N,\theta)_{\mathcal{H}} \geqslant_\ast (H,M,\varphi)_{\mathcal{H}}.
\] 
\end{lem}
\begin{proof}
  Suppose that $(\bar\Pj,\bar\Pj')$ gives $(\bar G,\bar N,\bar \theta)_{\hH}\geqslant_{\ast}(\bar H,\bar M,\bar \varphi)_{\hH}.$
  Notice that $\overline{G_{\theta}}=\bar G_{\bar\theta}$ and $\overline{H_{\varphi}}=\bar H_{\bar{\varphi}}.$
  Let $\Pj:G_{\theta}\ra\GL_{\theta(1)}(F)$ and $\Pj':H_{\varphi}\ra\GL_{\varphi(1)}(F)$ be inflations of $\bar\Pj$ and $\bar\Pj',$ respectively. 
 A routine verification shows that $(\Pj,\Pj')$ gives   $(G,N,\theta)_{\hH}\geqslant_{\ast}(H,M,\varphi)_{\hH}$ for the cases  $\ast\in\{\emptyset,c\}.$ 
  When $\ast=b,$ we need to check that
   \begin{equation}\label{equ:H-tri.lift.1}
    \Pj(\Cl_{\langle N,x \rangle}(x)^+)\text{ and }\Pj((\Cl_{\langle N,x \rangle}(x)\cap H)^+)\text{ are associated with a same scalar, } \forall x\in G_{\theta}.
   \end{equation}
  For any $x\in G_{\theta},$ there exists a positive integer $\lambda=|\Cl_{\langle N,x \rangle}(x)\cap Zx|$ such that \[\Pj(\Cl_{\langle N,x \rangle}(x)^+)=\lambda^{\ast}\bar \Pj(\Cl_{\langle \bar N,\bar x \rangle}(\bar x)^+)\] and \[\Pj((\Cl_{\langle N,x \rangle}(x)\cap H)^+)=\lambda^{\ast}\bar \Pj((\Cl_{\langle \bar N ,\bar x \rangle}(\bar x)\cap\bar H)^+),\]
  where  $\lambda^{\ast}$ is the image of $\lambda$ in $F.$
  Then (\ref{equ:H-tri.lift.1}) follows from the assumption that $\Pj(\Cl_{\langle \bar N,\bar x \rangle}(\bar x)^+)$ and $\Pj((\Cl_{\langle \bar N ,\bar x \rangle}(\bar x)\cap\bar H)^+)$ are associated with a same scalar.
\end{proof}

\begin{lem}[descending the order relations]\label{lem:H-tri.desc.}
Let $(G,N,\theta)_{\mathcal{H}}$ and $(H,M,\varphi)_{\mathcal{H}}$ be $\mathcal{H}$-triples satisfying $(G,N,\theta)_{\mathcal{H}} \geqslant_\ast (H,M,\varphi)_{\mathcal{H}}$ for $\ast \in \{\emptyset, c, b\}$. Consider a normal subgroup $Z \zg G$ contained in $\ker\theta \cap \ker\varphi$, and set $\bar{G} = G/Z$. Let $D$ be a defect group of $\bl(\bar{\varphi}).$ 
Assume $\overline{\C_G(N)}=\C_{\bar{G}}(\bar N)$ and $\C_{\bar G_{\bar\theta}}(D)\subseteq \bar H.$
Then 
\[
(\bar{G}, \bar{N}, \bar{\theta})_{\mathcal{H}} \geqslant_\ast (\bar{H}, \bar{M}, \bar{\varphi})_{\mathcal{H}}
\]
holds in both of the following cases:
\begin{enumerate}
    \item $Z$ is a $p'$-subgroup;
    \item $Z$ is a central $p$-subgroup of $G_\theta$.
\end{enumerate}
\end{lem}
\begin{proof}
  Suppose that  $(\Pj,\Pj')$ gives $(G,N,\theta)_{\hH}\geqslant_{\ast}(H,M,\varphi)_{\hH}.$
  Then $\Pj$ and $\Pj'$ descend to projective representations $\bar{\Pj}$ of $\bar{G}_{\bar \theta}$ and $\bar{\Pj}'$ of ${\bar H}_{\bar\varphi},$ respectively.
  A direct computation shows that $(\bar{\Pj},\bar{\Pj}')$ gives $(\bar G,\bar N,\bar \theta)_{\hH}\geqslant_{\ast}(\bar H,\bar M,\bar \varphi)_{\hH}$ when $\ast\in\{\emptyset,c\}.$
  
  Now assume $\ast=b.$
   Let $\nu$ be the isomorphism of character triples corresponding to $(\Pj,\Pj').$
   According to Theorem \ref{thm:2.two.equiv.}, it suffices to verify 
\[
 \bl(\overline{\nu_J(\chi)})^{\bar J}=\bl(\bar \chi)
\]
for all $J \in \mathcal{S}(G_\theta, N)$ and $\chi \in \IBr(J \which \theta)$. 
 This follows from \cite[Proposition 2.4]{NS14}, and the assumption that $\bl(\nu_J(\chi))^J=\bl(\chi).$  
\end{proof}

 The lemmas that follow are for proving Theorem \ref{thm:ibGAW.con.}, following \cite[Section 6]{FFZ}.
\begin{lem}\label{lem:H-tri1}
Let $(G,N,\theta)_{\hH}\geqslant_{\ast}(H,M,\varphi)_{\hH}$ be $\hH$-triples, where $\ast\in\{c,b\},$ and assume that $G$ is subgroup of a group $A.$ Then for any $a=(\sigma,g)\in\hH\times A,$ we have \[(G^a,N^a,\theta^a)_{\hH}\geqslant_\ast(H^a,M^a,\varphi^a).\]
\end{lem}
\begin{proof}
 The case $\ast = c$ was established in \cite[Lemma 6.1]{FFZ}. For $\ast = b$, we adapt the same approach as in \cite[Lemma 6.1]{FFZ}. 
Crucially, we note that the results of \cite[Lemmas 2.1, 2.3]{NSV20} remain valid for the partial order relation $\geqslant_b$ on $\mathcal{H}$-triples.
\end{proof}

We solve the problem at the end of the last section.

\begin{prop}\label{prop:conjB}
  Let $G$ be a finite group. 
  If Conjecture \ref{ConjB} holds for every block $B$ of $G,$ then it holds for $G.$
\end{prop}
\begin{proof}
Through appropriate adjustments of the bijection $\Omega:\IBr(G)\ra \W^{\circ}(G)/G,$ we may assume $\Omega$ is $\hH\times \Aut(G)$-equivariant.
Then we  apply Lemma \ref{lem:H-tri1}.
\end{proof}

Let $N$ be a non-trivial finite group, and $\bS_m$ the symmetric group  of degree $m \geq 1.$
The wreath product $\Aut(N) \wr \bS_m$ embeds naturally as a subgroup of $\Aut(N^m)$, where $N^m = N \times \cdots \times N$ ($m$ factors) denotes the direct product.
See the discussion preceding \cite[Lemma 6.2]{FFZ} for details.

\begin{lem}\label{lem:H-tri2}
Let $(G,N,\theta)_{\mathcal{H}} \geqslant_\ast (H,M,\varphi)_{\mathcal{H}}$ ($\ast \in \{c,b\}$) be $\hH$-triples. 
Let $\sigma_1,\ldots,\sigma_k \in \mathcal{H}$ and $l \geq 1$ an integer. 
Denote $m=kl$ and set $\widetilde{\theta}= (\theta^{\sigma_1})^{l}\times\cdots \times(\theta^{\sigma_k})^{l} \in\IBr(N^m)$, $\widetilde{\varphi}=(\varphi^{\sigma_1})^{l}\times\cdots
\times(\varphi^{\sigma_k})^{l}\in\IBr(M^m)$. Assume that the $\theta^{\sigma_i},1\leq i\leq k,$ are pairwise not $G$-conjugate, then we have 
\[((G\wr \bS_m)_{\widetilde{\theta}^{\hH}},N^m,\widetilde{\theta})_{\mathcal{H}} \geqslant_\ast 
((H\wr \bS_m)_{\widetilde{\varphi}^{\hH}},M^m,
\widetilde{\varphi})_{\mathcal{H}}.\]
\end{lem}
\begin{proof}
The case $\ast = c$ was settled in \cite[Lemma 6.2]{FFZ}. Now consider $\ast = b$. 
Following the proof of \cite[Lemma 6.2]{FFZ} while replacing $\eta$ and $\eta'$ with $\theta$ and $\varphi$ respectively, we have established that
\begin{equation}\label{equ:H-tri4}
    \big((G \wr \bS_m)_{\widetilde{\theta}^{\hH}}, N^m, \widetilde{\theta}\big)_{\hH} \geqslant_c \big((H \wr \bS_m)_{\widetilde{\varphi}^{\hH}}, M^m, \widetilde{\varphi}\big)_{\hH}.
\end{equation}
Assume that \eqref{equ:H-tri4} is given by  $(\widetilde{\Pj}, \widetilde{\Pj}')$. To complete the proof, we must show that this  pair gives
\begin{equation}\label{equ:H-tri5}
    \big((G \wr \bS_m)_{\widetilde{\theta}}, N^m, \widetilde{\theta}\big) \geqslant_b \big((H \wr \bS_m)_{\widetilde{\varphi}}, M^m, \widetilde{\varphi}\big).
\end{equation}
Following the construction in \cite[Theorem 2.7]{NSV20}, we recall that
$\widetilde{\Pj}=\widehat{\Pj}^{\sigma_1}\otimes \cdots\otimes\widehat{\Pj}^{\sigma_k}\text{ and }
  \widetilde{\Pj}'=\widehat{\Pj}'^{\sigma_1}\otimes \cdots\otimes\widehat{\Pj}'^{\sigma_k},$
where $(\widehat{\Pj}, \widehat{\Pj}')$ gives
\begin{equation}\label{equ:H-tri6}
    \big(G_{\theta} \wr \bS_l, N^l, \theta^l\big) \geqslant_c \big(H_{\varphi} \wr \bS_l, M^l, \varphi^l\big).
\end{equation}
Moreover, as demonstrated in \cite[Lemma 3.10]{MRR23}, the pair $(\widehat{\Pj}, \widehat{\Pj})$ provides a block isomorphism of character triples. 
The relation \eqref{equ:H-tri5} follows immediately from Definition \ref{de:char-tri} by direct computation.
\end{proof}

\begin{lem}\label{lem:H-tri3}
  Let $m\geq 1$ be an integer. 
  For each $1 \leq i \leq m$, let $(G_i,N_i,\theta_i)_{\hH}$ and $(H_i,M_i,\varphi_i)_{\hH}$ be $\hH$-triples such that $(G_i,N_i, \theta_i)_{\hH} \geqslant_\ast (H_i,M_i, \varphi_i)_{\hH}$, where $\ast \in \{c, b\}$. 
  Let $X=X_1\times\cdots\times X_m, X\in\{G,N,H,M,\theta,\varphi\}.$ Then \[(G_{\theta^{\hH}},N,\theta)_{\hH}\geqslant_{\ast} (H_{\varphi^{\hH}},M,\varphi)_{\hH}.\]
\end{lem}
\begin{proof}
  The case of $\ast=c$ has been settled by \cite[Lemma 6.3]{FFZ}.
  Let $\ast=b.$
  We only need to verify the extra condition of the block isomorphism of character triples in the proof of \cite[Lemma 6.3]{FFZ}.
  And this is easily seen, see also \cite[Lemma 3.9]{MRR23}.
\end{proof}

\section{Centralization of $\hH$-triples}
This section is devoted to the proof of Theorem~\ref{thm:centrali.}, 
a key ingredient in our main reduction.
We first prove two preliminary lemmas.

\begin{lem}\label{lem:cen_pre}
Let $N \unlhd G$ and $Q/N \unlhd G/N$ a $p$-group. If $B$ is a block of $G$ covering a $G$-invariant block $b$ of $N$, then $Q \subseteq DN$ for any defect group $D$ of $B$.
\end{lem}
\begin{proof}
Let $\hat{b}$ be the unique block of $Q$ covering $b$. By Fong's theorem \cite[Chapter V, Theorem 5.16]{NagTsu}, there exists a defect group $D_1$ of $\hat{b}$ with $D_1N = Q$. Since $B$ covers $\hat{b}$, Fong's theorem guarantees a defect group $D$ of $B$ containing $D_1$. Therefore, $Q = D_1N \subseteq DN$.
\end{proof}

\begin{lem}\label{lem:vert.}
Let $N \zg G$ be finite groups and write $\bar G = G/N$. The following two results on vertices of irreducible Brauer characters hold.
\begin{enumerate}
  \item 
  Let $\theta \in \IBr(\bar G)$ and $\hat\theta \in \IBr(G)$ be its inflation through the canonical epimorphism.
  If $V$ is a vertex of $\hat\theta$, then $\bar V$ is a vertex of $\theta$, and $V \cap N$ is a Sylow $p$-subgroup of $N$.
  
  \item Let $\zeta \in \IBr(N)$ be $G$-invariant, and let $\varphi \in \IBr(G)$ be such that $\zeta \varphi_N \in \IBr(N)$.
      Then for any $\chi \in \IBr(G \which \zeta)$ and any vertex $V$ of $\chi$, there exists a vertex $V_1$ of $\chi \varphi$ such that $VN = V_1N$.
\end{enumerate}
\end{lem}
\begin{proof}
For (1), let $X: \bar G \to \GL_{\theta(1)}(F)$ be a group representation affording the character $\theta$, and let $\hat X: G \to \GL_{\theta(1)}(F)$ be its natural inflation.
Regard $\GL_{\theta(1)}(F)$ as the unit group of $\M_{\theta(1)}(F)$.
Then the group homomorphisms $X$ and $\hat X$ make $\M_{\theta(1)}(F)$ a $\bar G$-algebra and a $G$-algebra, denoted by $A$ and $B$, respectively.
By the definition of vertices and \cite[Chapter IV, Theorem 2.2]{NagTsu}, a vertex of $\hat\theta$ is a minimal subgroup $V$ of $G$ such that $1_B \in B_V^G$ (where $B_V^G$ is the image of the trace map $\Tr_V^G: B^V \to B^G$; see the statements preceding \cite[Lemma 2.3]{Fu24} for details).

If $V$ is a vertex of $\hat\theta$, then $1_B \in B_V^G$. Thus, $1_B \in B_{VN}^G$ and $1_A \in A_{\bar V}^{\bar G}$.
Let $V_1$ be a minimal subgroup of $\bar G$ such that $1_A \in A_{V_1}^{\bar G}$ and $V_1 \subseteq \bar V$.
Then $V_1$ is a vertex of $\theta$. We claim that $V_1 = \bar V$. Otherwise, $V_1 \lneqq \bar V$.
Let $V' \in \mS(G, N)$ be such that $V'/N = V_1$.
Since $1_A \in A_{V_1}^{\bar G}$, we have $1_B \in B_{V'}^G$.
By \cite[Lemma 2.3]{Fu24}, we have $V^g \subseteq V'$ for some $g \in G$, a contradiction.
The fact that $V \cap N$ is a Sylow $p$-subgroup of $N$ follows from \cite[Chapter IV, Lemma 3.4(ii)]{NagTsu}, since $\hat\theta$ lies over $1_N$ (the trivial character of $N$), and the vertices of $1_N$ are the Sylow $p$-subgroups of $N$.

For (2), let $\Pj$ be a projective representation of $G$ associated with $\zeta$.
Let $\PE$ be a group representation of $G$ affording $\varphi$, and let $\Pj' = \Pj \otimes \PE$.
Then $\Pj'$ is a projective representation of $G$ associated with $\zeta \varphi_N$.
The pair $(\Pj, \Pj')$ induces a bijection 
\[
f: \IBr(G \which \zeta) \to \IBr(G \which \zeta \varphi_N), \quad \chi \mapsto \chi \varphi.
\]
By \cite[Proposition 2.4]{Fu24}, the bijection $f$ satisfies the required condition.
\end{proof}

As established in \cite[Section 3]{Fu24}, every $\hH$-triple has an associated cohomology class. We now summarize this construction.
Let $^*: \mO \to F$ denote the natural ring homomorphism and $\Fp$ be the prime subfield of $F$ with $p$ elements. Given an $\hH$-triple $(G,N,\theta)_{\hH}$, define $E := \Fp[\theta]$ as the field extension of $\Fp$ generated by $\{\theta(x)^* \mid x \in N_{p'}\}$, $A := \M_{\theta(1)}(E)$ the matrix algebra over $E$, and $X: N \to \GL_{\theta(1)}(E) = A^\times$ a representation affording $\theta$ and realized over $E$.
Let $s = [E:\Fp]$ and fix an $\Fp$-algebra embedding $\iota: A \hookrightarrow \M_{s\theta(1)}(\Fp)$ preserving identities. A \emph{projective representation of $G$ associated with $\theta$} is a map $Y: G \to \GL_{s\theta(1)}(\Fp)$ satisfying:
\begin{enumerate}
    \item $Y(n) = X(n)$ for all $n \in N$;
    \item $Y(gn) = Y(g)Y(n)$ and $Y(ng) = Y(n)Y(g)$ for all $n \in N$, $g \in G$.
\end{enumerate}
The map $\alpha: G \times G \to E^\times$ defined by $Y(g)Y(h) = Y(gh)\alpha(g,h)$ is a factor set in $\Z^2(G,E^\times)$,  where $g\in G$ acts on $E$ via automorphism $\sigma$ such that $(\sigma, g)\in(\hH\times G)_{\theta}.$ By \cite[Lemma 3.1]{Fu24}, we have $\C_G(E) = G_\theta$.
Furthermore, 
 $\alpha$ is constant on $N \times N$-cosets, and
 $\alpha(1,g) = \alpha(g,1) = 1$ for all $g \in G$.
We say $\alpha$ is \emph{a factor set associated with $(G,N,\theta)_{\hH}$}.

\begin{thm}\label{thm:centrali.}
  Let $(A,N,\theta)_{\hH}$ be an $\hH$-triple. There exists an $\hH$-triple $(K,Z,\lambda)_{\hH}$ satisfying the following conditions.
  \begin{enumerate}
    \item There exists a group isomorphism $\epsilon:A/N\ra K/Z$ and we identify $\bar{A}=A/N$ with $K/Z.$ The subgroup $Z$ is $p'$-central in $K_{\lambda}$ and $\lambda$ is a faithful linear character of $Z.$ For any $J\in\mS(A,N),$ define $J^{\bullet}\leq K$ via $J^{\bullet}/Z=\epsilon(J/N).$
    \item $(\hH\times \bar A)_{\theta}=(\hH\times \bar A)_{\lambda}.$ 
    \item For each $J\in\mS(A_{\theta},N),$ there exists a bijection $$\nu_J:\IBr(J\which \theta)\ra \IBr(J^{\bullet}\which\lambda)$$ 
        with the following properties:
\begin{enumerate}
  \item   $\nu_J(\chi)^a=\nu_{J^a}(\chi^a)$ for all $\chi\in\IBr(J\which\theta)$  and $a\in(\hH\times \bar A)_{\theta}$;
  \item   For any vertex $V$ of $\chi\in\IBr(J\which\theta),$ there exists a vertex $V_1$ of $\nu_{J}(\chi)$ satisfying $\bar V=V_1Z/Z$.
\end{enumerate}
We denote $\nu_J(\chi)$ by $\chi^{\bullet}.$
 
    \item     Let $N\subseteq Q\subseteq J\subseteq G$ be subgroups of $A$ with $Q/N$ a $p$-group and $J\subseteq A_{\theta}.$ 
        Let $\chi\in\IBr(J\which\theta)$ and $\varphi\in\IBr(\N_J(Q)\which \theta).$  
        Let $\ast\in\{\emptyset,c,b\}.$ 
        We write $H=\N_G(Q)$ and $M=\N_J(Q).$
        If $(G^{\bullet},J^{\bullet},\chi^{\bullet})_{\hH}$ and $(H^{\bullet},M^{\bullet},\varphi^{\bullet})_{\hH}$ are $\hH$-triples and $$(G^{\bullet},J^{\bullet},\chi^{\bullet})_{\hH}
        \geqslant_{\ast}
         (H^{\bullet},M^{\bullet},\varphi^{\bullet})_{\hH},$$
          then 
          $$(G,J,\chi)_{\hH}
        \geqslant_{\ast}
         (H,M,\varphi)_{\hH}.$$
  \end{enumerate}
  \end{thm}
\begin{proof}
Let \( E = \Fp[\theta] \), and let \( \alpha \in \Z^2(A, E^{\times}) \) be a factor set associated with the triple \( (A, N, \theta)_{\hH} \).  
The set  
\(  
\widehat{A} = \{(g, z) \which g \in A, \, z \in E^{\times}\}  
\) 
forms a finite group under the multiplication  
\[  
(g, z_1)(h, z_2) = \big(gh, \, \alpha(g, h) z_1^{h} z_2\big), \quad \forall (g, z_1), (h, z_2) \in \widehat{A}.  
\]  
Let \( Z_0 = \{(1, z) \which z \in E^{\times}\} \).  
Then \( Z_0 \) is a normal \( p' \)-subgroup of \( \widehat{A} \), and \( \widehat{A}/Z_0 \cong A \).  
We identify \( A \) with \( \widehat{A}/Z_0 \) via the natural isomorphism.  
For any \( J \in \mS(A, N) \), we denote by \( \widehat{J} \leq \widehat{A} \) the subgroup satisfying \( J = \widehat{J}/Z_0 \).  
Note that \( Z_0 \) is a central subgroup of \( \widehat{A_{\theta}} \),  
and \( \widehat{N} = Z_0 \times N_0 \), where \( N_0 = \{(n, 1) \which n \in N\} \) is normal in \( \widehat{A} \) and isomorphic to \( N \).

Let \(\theta_0 \in \IBr(N_0)\) be the character corresponding to \(\theta \in \IBr(N)\) via the natural isomorphism.  
Let \(\Pj\) be a projective representation of \(A_\theta\) associated with \(\theta\), with factor set \(\alpha_{A_\theta \times A_\theta}\).  
By \cite[Lemma 3.7]{Fu24}, for any \(a = (\sigma, g) \in (\hH \times A)_\theta\), the function \(\mu_a \colon A_\theta \to F^\times\) determined by \(\Pj^a \sim \mu_a \Pj\) satisfies  
\[
\mu_a(h) = \Big( \alpha(g, h g^{-1})^{-1} \alpha(h, g^{-1})^{-1} \alpha(g, g^{-1}) \Big)^\sigma, \quad \forall h \in A_\theta.
\]
There exists a natural extension \(\widetilde{\theta} \in \IBr(\widehat{A_\theta})\) of \(\theta_0\), afforded by the group representation  
\[
\widehat{\Pj} \colon \widehat{A_\theta} \to \GL_{\theta(1)}(F), \quad (h, z) \mapsto z \Pj(h).
\]  
Let \(\hat{\theta} \in \IBr(\widehat{N})\) be the inflation of \(\theta\) via the natural epimorphism \(\widehat{N} \to N\).  
One verifies directly that
\[
(\sigma, (g, z)) \in (\hH \times \widehat{A})_{\hat{\theta}} \iff (\sigma, g) \in (\hH \times A)_\theta.
\]  
We now claim that \(\widetilde{\theta}\) is \((\hH \times \widehat{A})_{\hat{\theta}}\)-invariant.  
Let \(\hat{a} = (\sigma, (g, z_1)) \in (\hH \times \widehat{A})_{\hat{\theta}}\) and $a=(\sigma, g)\in(\hH\times A)_{\theta}$.  
For any \((h, z) \in \widehat{A_\theta}\), we compute  
\begin{align*}
\widehat{\Pj}^{\hat{a}}(h, z) 
&= \widehat{\Pj}\big( (g, 1)(h, z)(g, 1)^{-1} \big)^\sigma \\
&= \widehat{\Pj}\Big( ghg^{-1}, \, \alpha(g, h g^{-1}) \alpha(h, g^{-1}) \alpha(g, g^{-1})^{-1} z^{g^{-1}} \Big)^\sigma \\
&= z \mu_a(h)^{-1} \Pj^a(h) \\
&\sim z \Pj(h) = \widehat{\Pj}(h, z),
\end{align*}  
as claimed.

Let $K = \widehat{A}/N_0$ and $Z = Z_0N_0/N_0$. Let $\upsilon \in \IBr(Z_0)$ be the linear character afforded by the representation
\[
Z_0 \to F^{\times}, \quad (1,z) \mapsto z.
\]
Define $\lambda_1 = \upsilon^{-1} \times \mathbf{1}_{N_0} \in \IBr(Z_0 \times N_0)$, and let $\lambda \in \IBr(Z)$ be its deflation through the natural epimorphism $Z_0 \times N_0 \to Z$.
Then $Z$ is a $p'$-central subgroup of $\widehat{A_{\theta}}/N_0$, and $\lambda$ is a faithful linear character of $Z$.
There is a well-defined isomorphism
\[
\epsilon \colon A/N \to K/Z, \quad gN \mapsto (g,1)Z_0N_0,
\]
and for any $J \in \mS(A,N)$, we have $J^\bullet = \widehat{J}/N_0$.

The representation affording $\lambda$ is
\[
X_\lambda \colon Z \to F^{\times}, \quad (1,z)N_0 \mapsto z^{-1}.
\]
For $(\sigma, (g,z')N_0) \in \hH \times K$, we compute:
\begin{align*}
X_\lambda^{(\sigma,(g,z')N_0)}\big((1,z)N_0\big) 
&= X_\lambda\left( (g,1)(1,z)(g,1)^{-1}N_0 \right)^\sigma \\
&= X_\lambda\left( (1,z^{g^{-1}})N_0 \right)^\sigma \\
&= \left(z^{g^{-1}}\right)^{-\sigma}.
\end{align*}
Since an element $(\sigma,g) \in \hH \times A$ stabilizes $\theta$ if and only if the restriction $\sigma_E$ coincides with the action of $g$ on $E$,
it follows that $(\sigma,(g,z')N_0)$ stabilizes $\lambda$ precisely when $(\sigma,g)$ stabilizes $\theta$.
This implies that $K_\lambda = \widehat{A_\theta}/N_0$, and completes the proof of  (1) and (2).

For any $J \in \mS(A_\theta, N)$, there exist natural bijections induced by inflation:
\begin{align*}
\nu_{J,1} &: \IBr(J \which \theta) \to \IBr(\widehat{J} \which \hat\theta), \\
\nu_{J,2} &: \IBr(J^\bullet \which \lambda) \to \IBr(\widehat{J} \which \lambda_1).
\end{align*}
Since $\lambda_1\widetilde{\theta}_{\widehat{N}} = \hat\theta$, Corollary 8.19 in \cite{Nav:char.} yields a bijection
\[
\nu_{J,3} : \IBr(\widehat{J} \which \lambda_1) \to \IBr(\widehat{J} \which \hat\theta), \quad \zeta \mapsto \zeta\widetilde{\theta}_{\widehat{J}}.
\]
By part (2), we have $(\hH \times \widehat{A})_{\hat\theta} = (\hH \times \widehat{A})_{\lambda_1}$. For any $a = (\sigma, (g, z)) \in (\hH \times \widehat{A})_{\hat\theta}$, the $(\hH \times \widehat{A})_{\hat\theta}$-invariance of $\widetilde{\theta}$ implies
\[
\nu_{J,3}(\zeta)^a = (\zeta\widetilde{\theta}_{\widehat{J}})^a = \zeta^a \cdot \widetilde{\theta}_{\widehat{J}^a} = \nu_{J^a,3}(\zeta^a),
\]
where $J^a = J^g$.
We  define the composite bijection
\[
\nu_J : \IBr(J \which \theta) \to \IBr(J^\bullet \which \lambda), \quad \nu_J = \nu_{J,2}^{-1}\nu_{J,3}^{-1} \nu_{J,1}.
\]
This bijection $\nu_J$ is $(\hH\times \bar A)_{\theta}$-equivariant by the above arguments.
Let $V \leq J$ be a vertex of $\chi \in \IBr(J \which \theta)$. 
Then by Lemma \ref{lem:vert.}, there exists a vertex $V_1 \leq \widehat J$ of $\nu_{J,1}(\chi)$ such that $V = V_1 Z_0 / Z_0$.
Also by Lemma \ref{lem:vert.}, there exists a vertex $V_2$ of $\nu_{J,3}^{-1}\big( \nu_{J,1}(\chi)\big)$ such that $V_1 \widehat{N} = V_2 \widehat{N}$.
Since $V_2 N_0 / N_0$ is a vertex of $\nu_J(\chi)$, statement (3.b) follows by a direct group-theoretic computation.

Let $\zeta \in \IBr(\widehat{J})$ and $\delta \in \IBr(\widehat{M})$ be the inflations of $\chi^{\bullet}$ and $\varphi^{\bullet}$ through the corresponding epimorphisms, respectively.
Then $\hat\chi := \zeta \widetilde{\theta}_{\widehat{J}}$ and $\hat\varphi := \delta \widetilde{\theta}_{\widehat{M}}$ are the inflations of $\chi$ and $\varphi$ through the corresponding epimorphisms, respectively.
Since by assumption
\begin{equation}\label{equ:cen.0}
  (G^{\bullet}, J^{\bullet}, \chi^{\bullet})_{\hH}
        \geqslant_{\ast}
         (H^{\bullet}, M^{\bullet}, \varphi^{\bullet})_{\hH},
\end{equation}
by Lemma \ref{lem:H-tri.lift}, the above relation naturally lifts to
\begin{equation}\label{equ:cen.1}
 (\widehat{G}, \widehat{J}, \zeta)_{\hH} \geqslant_{\ast} (\widehat{H}, \widehat{M}, \delta)_{\hH}.
\end{equation}

Suppose that \eqref{equ:cen.1} is given by $(\PQ, \PE)$, where $\PQ$ and $\PE$ are projective representations of $\widehat{G}_{\zeta}$ and $\widehat{H}_{\delta}$, respectively.
Note that $\widehat{G}_{\zeta}$ and $\widehat{H}_{\delta}$ are contained in $\widehat{A_{\theta}}$.
We now prove that $(\PQ \otimes \widehat{\Pj}_{\widehat{G}_{\zeta}}, \PE \otimes \widehat{\Pj}_{\widehat{H}_{\delta}})$ gives
\begin{equation}\label{equ:cen.2}
  (\widehat{G}, \widehat{J}, \zeta \widetilde{\theta}_{\widehat{J}})_{\hH} \geqslant_{\ast} (\widehat{H}, \widehat{M}, \delta \widetilde{\theta}_{\widehat{M}})_{\hH},
\end{equation}
where $\widehat{\Pj}$ is the group representation affording $\widetilde{\theta}$.

First, note that
\[
(\hH \times \widehat{H})_{\zeta} = (\hH \times \widehat{H})_{\delta} = (\hH \times \widehat{H})_{\hat\chi} = (\hH \times \widehat{H})_{\hat\varphi} \subseteq (\hH \times \widehat{A})_{\hat{\theta}} = (\hH \times \widehat{A})_{\lambda_1}.
\]
For any $a \in (\hH \times \widehat{H})_{\zeta}$, since $\widetilde{\theta}^a = \widetilde{\theta}$, we have $\widehat{\Pj}^a \sim \widehat{\Pj}$.
Let $\mu_a'\colon \widehat{G}_{\zeta} \to F^{\times}$ and $\mu_a''\colon \widehat{H}_{\delta} \to F^{\times}$ be the functions determined by $\PQ^a \sim \mu_a' \PQ$ and $\PE^a \sim \mu_a'' \PE$, respectively.
Note that $\mu_a'' = (\mu_a')_{\widehat{H}_{\delta}}$.
Then
\[
(\PQ \otimes \widehat{\Pj}_{\widehat{G}_{\zeta}})^a = \PQ^a \otimes (\widehat{\Pj}^a)_{\widehat{G}_{\zeta}} \sim \mu_a' \PQ \otimes \widehat{\Pj}_{\widehat{G}_{\zeta}}.
\]
Similarly, we have
\[
(\PE \otimes \widehat{\Pj}_{\widehat{H}_{\delta}})^a \sim \mu_a'' \PE \otimes \widehat{\Pj}_{\widehat{H}_{\delta}}.
\]
For the case when $\ast \in \{\emptyset, c\}$, the remaining conditions follow by direct computation.

Now we discuss the case $\ast = b$.
Since $\widehat{Q}/\widehat{N}$ is a normal $p$-subgroup of $\widehat{M}/\widehat{N}$, and $\hat\varphi$ lies over the $\widehat{M}$-invariant character $\hat\theta$, by Lemma \ref{lem:cen_pre}, a defect group $D$ of $\bl(\hat\varphi)$ satisfies $\widehat{Q} \subseteq D \widehat{N}$.
Hence, $\C_{\widehat{G}}(D) \subseteq \N_{\widehat{G}}(\widehat{Q}) = \widehat{\N_G(Q)} = \widehat{H}$.
For $L \in \mS(G_{\chi}, J)$, let
\[
\kappa_{\widehat{L}} \colon \IBr(\widehat{L} \which \zeta) \to \IBr(\widehat{L} \cap \widehat{H} \which \delta)
\]
be the bijection induced by $(\PQ_{\widehat{L}}, \PE_{\widehat{L} \cap \widehat{H}})$.
Then the bijection $\IBr(\widehat{L} \which \hat\theta) \to \IBr(\widehat{L} \cap \widehat{H} \which \hat\varphi)$ induced by $(\PQ_{\widehat{L}} \otimes \widehat{\Pj}_{\widehat{L}}, \PE_{\widehat{L} \cap \widehat{H}} \otimes \widehat{\Pj}_{\widehat{L} \cap \widehat{H}})$ sends $\rho \widetilde{\theta}_{\widehat{L}}$ to $\kappa_{\widehat{L}}(\rho) \widetilde{\theta}_{\widehat{L} \cap \widehat{H}}$ for all $\rho \in \IBr(\widehat{L} \which \zeta)$.
We need to prove that
\begin{equation}\label{equ:cen.3}
  \bl\big( \kappa_{\widehat{L}}(\rho) \widetilde{\theta}_{\widehat{L} \cap \widehat{H}} \big)^{\widehat{L}} = \bl(\rho \widetilde{\theta}_{\widehat{L}}) \quad \text{for all } \rho \in \IBr(\widehat{L} \which \zeta).
\end{equation}

Let $\rho\in  \IBr(\widehat{L} \which \zeta).$
For any fixed $k \in \widehat{L}$, we define an equivalence relation on the conjugacy class $\Cl_{\widehat{L}}(k)$ such that $x, x' \in \Cl_{\widehat{L}}(k)$ are equivalent if and only if they belong to the same $N_0$-coset. Let
\[
\Cl_{\widehat{L}}(k) = \mathcal{C}_1 \sqcup \cdots \sqcup \mathcal{C}_s \sqcup \mathcal{C}_{s+1} \sqcup \cdots \sqcup \mathcal{C}_t
\]
be the corresponding partition, with $\mathcal{C}_i \subseteq \widehat{L} \cap \widehat{H}$ for $1 \leq i \leq s$, and $\mathcal{C}_i \subseteq \widehat{L} \setminus \widehat{H}$ for $s+1 \leq i \leq t$.
Since $N_0$ is normal in $\widehat{L}$, we have $|\mathcal{C}_1| = \cdots = |\mathcal{C}_t| =: r$.
Choose $x_i \in \mathcal{C}_i$ for each $i$; then $\{x_i N_0 \which 1 \leq i \leq t\}$ is the conjugacy class $\Cl_{L^{\bullet}}(k N_0)$.
Let
\[
\rho_1 \in \IBr(L^{\bullet}) \quad \text{and} \quad \kappa_{\widehat{L}}(\rho)_1 \in \IBr((L \cap H)^{\bullet})
\]
be the deflations of $\rho$ and $\kappa_{\widehat{L}}(\rho)$ under the natural epimorphisms, respectively (note that $(L \cap H)^{\bullet} = L^{\bullet} \cap H^{\bullet} = (\widehat{L} \cap \widehat{H}) / N_0$).
Then by assumption \eqref{equ:cen.0}, we have
\begin{equation}\label{equ:cen.4}
  \bl(\kappa_{\widehat{L}}(\rho)_1)^{L^{\bullet}} = \bl(\rho_1).
\end{equation}
Let $X_{\rho}$ and $X_{\rho}'$ be the representations of $\widehat{L}$ and $\widehat{L} \cap \widehat{H}$ affording $\rho$ and $\kappa_{\widehat{L}}(\rho)$, respectively.
Then \eqref{equ:cen.4} implies that
\begin{equation}\label{equ:cen.5}
  \sum_{i=1}^{t} X_{\rho}(x_i) \quad \text{and} \quad \sum_{i=1}^{s} X_{\rho}'(x_i) \quad \text{are associated with the same scalar}.
\end{equation}
Note that $X_{\rho} \otimes \widehat{\Pj}_{\widehat{L}}$ and $X_{\rho}' \otimes \widehat{\Pj}_{\widehat{L} \cap \widehat{H}}$ afford $\rho \widetilde{\theta}_{\widehat{L}}$ and $\kappa_{\widehat{L}}(\rho) \widetilde{\theta}_{\widehat{L} \cap \widehat{H}}$, respectively.
Since $\widehat{\Pj}_{N_0}$ is irreducible and each $\mathcal{C}_i$ ($1 \leq i \leq t$) is a union of $N_0$-conjugacy classes, each $\widehat{\Pj}(\mathcal{C}_i^{+})$ is a scalar matrix.
As any $\mathcal{C}_i$ ($1 \leq i \leq t$) is $\widehat{L}$-conjugate to $\mathcal{C}_1$, the matrices $\widehat{\Pj}(\mathcal{C}_i^+)$ for $1 \leq i \leq t$ are all associated with the same scalar, say $\epsilon$.
We compute:
\begin{align*}
 &\big( X_{\rho} \otimes \widehat{\Pj}_{\widehat{L}} \big) (\Cl_{\widehat{L}}(k)^+) 
 = \sum_{x \in \Cl_{\widehat{L}}(k)} X_{\rho}(x) \otimes \widehat{\Pj}(x) 
 = \sum_{i=1}^{t} \left( \sum_{x \in \mathcal{C}_i} X_{\rho}(x) \otimes \widehat{\Pj}(x) \right) \\
 &= \sum_{i=1}^{t} \left( X_{\rho}(x_i) \otimes \widehat{\Pj}(\mathcal{C}_i^+) \right)
 = \epsilon \left( \sum_{i=1}^{t} X_{\rho}(x_i) \right) \otimes \I_{\theta(1)}.
\end{align*}
Similarly, we compute:
\[
\big( X_{\rho}' \otimes \widehat{\Pj}_{\widehat{L} \cap \widehat{H}} \big) ((\Cl_{\widehat{L}}(k) \cap \widehat{H})^+) 
= \epsilon \left( \sum_{i=1}^{s} X_{\rho}'(x_i) \right) \otimes \I_{\theta(1)}.
\]
Combining this with \eqref{equ:cen.5}, we conclude that
\[
\big( X_{\rho} \otimes \widehat{\Pj}_{\widehat{L}} \big) (\Cl_{\widehat{L}}(k)^+) \quad \text{and} \quad \big( X_{\rho}' \otimes \widehat{\Pj}_{\widehat{L} \cap \widehat{H}} \big) ((\Cl_{\widehat{L}}(k) \cap \widehat{H})^+)
\]
are associated with the same scalar.
By varying $k \in \widehat{L}$, we prove \eqref{equ:cen.3}, and thus \eqref{equ:cen.2}.

It is straightforward to compute that $\C_{\widehat{G}}(\widehat{J}) Z_0 / Z_0 = \C_G(J)$ (or apply \cite[Theorem 4.1 (d)]{NS14}).
By Lemma \ref{lem:cen_pre}, a defect group $D_1$ of $\bl(\varphi)$ satisfies $Q \subseteq D_1N$.
Hence, $\C_{G}(D_1) \subseteq \N_{G}(Q)  = H$.
Thus, we can apply Lemma \ref{lem:H-tri.desc.} to \eqref{equ:cen.2} and obtain
\[
(G, J, \chi)_{\hH} \geqslant_{\ast} (H, M, \varphi)_{\hH}.
\]
This completes the proof of the theorem.
\end{proof}

\section{Induction on $\hH$-triples}
Suppose that $(G, N, \theta)$ and $(H, M, \varphi)$ are character triples such that $G = NH$, $M = N \cap H$, and $\theta = \varphi^N$.
Let $\Pj$ be a projective representation of $H$ associated with $\varphi$, and let $\{n_1=1, \dots, n_s\}$ be a complete set of representatives of right $M$-cosets in $N$.
Then $\{n_1, \dots, n_s\}$ is also a complete set of representatives of right $H$-cosets in $G$.
Note that $\theta(1) = \varphi(1) s$.
We define the function
\[
\Ind^G_{H,N}(\Pj) \colon G \to \GL_{\theta(1)}(F), \quad x \mapsto 
\begin{pmatrix}
   \Pj(n_1 x n_1^{-1}) & \cdots & \Pj(n_1 x n_s^{-1}) \\
   \vdots & & \vdots \\
   \Pj(n_s x n_1^{-1}) & \cdots & \Pj(n_s x n_s^{-1})
\end{pmatrix},
\]
where we define $\Pj(n_i x n_j^{-1})$ to be the zero matrix if $n_i x n_j^{-1} \notin H$ for $i, j \in \{1, \dots, s\}$. The notation $\Ind^G_{H,N}(\Pj)$ follows \cite[Remark 3.1]{Fe25}, and the idea originates from the proof of \cite[Theorem 3.14]{NS14}.
As noted in the proof of \cite[Theorem 3.14]{NS14}, the function $\widehat{\Pj}$ is a projective representation of $G$ associated with $\theta$, and the factor sets of $\Pj$ and $\widehat{\Pj}$ coincide under the isomorphism $G/N \cong H/M$.

Let $\bS_s$ be the symmetric group of degree $s$, and assume that the elements in $\bS_s$ are composed from the left (i.e., $w v(i) = v(w(i))$ for any $w, v \in \bS_s$).
For any $x \in G$, let $w_x \in \bS_s$ be defined by $H n_i x = H n_{w_x(i)}$.
It is straightforward to verify that $w_x w_y = w_{xy}$ for any $x, y \in G$.
For any $w \in \bS_s$, define $T_w := T_{w, \varphi(1)} \in \GL_{\theta(1)}(F)$ to be the permutation matrix such that for any $i, j \in \{1, \dots, s\}$, the $(i,j)$-th block of $T_w$ is the identity matrix of degree $\varphi(1)$ if $j = w(i)$, and the zero matrix of degree $\varphi(1)$ otherwise.
Note that $T_w T_v = T_{w v}$ for any $w, v \in \bS_s$.

Let $\widehat{\Pj} = \Ind^G_{H,N}(\Pj)$. 
Then for any $x \in G$, the matrix $\widehat{\Pj}(x)$ can be expressed as
\[
\widehat{\Pj}(x) = \diag\big( \Pj(n_1 x n_{w_x(1)}^{-1}), \dots, \Pj(n_s x n_{w_x(s)}^{-1}) \big) T_{w_x}.
\]
Note that for any block diagonal matrix $\diag(A_1, \dots, A_s) \in \M_{\theta(1)}(F)$ with $A_1, \dots, A_s \in \M_{\varphi(1)}(F)$, we have
\[
T_{w}^{-1} \diag(A_1, \dots, A_s) T_w = \diag(A_{w^{-1}(1)}, \dots, A_{w^{-1}(s)}).
\]

\begin{thm}\label{thm:H-tri(ind)}
Suppose that $(G, N, \theta)_{\hH}$ and $(H, M, \varphi)_{\hH}$ are $\hH$-triples such that $G = NH$, $N \cap H = M$, and $(\hH \times H)_{\theta} = (\hH \times H)_{\varphi}$.
Let $K \leq N$, $Z \leq M$ be such that $Z = K \cap H$.
Let $\chi \in \IBr(K)$, $\lambda \in \IBr(Z)$ be such that $\chi^N = \theta$, $\lambda^M = \varphi$, and $N_{\chi} = K$, $M_{\lambda} = Z$.
Let $\ast \in \{\emptyset, c, b\}$.
Suppose that
\[
(G_{\chi^{\hH}}, K, \chi)_{\hH} \geqslant_{\ast} (H_{\lambda^{\hH}}, Z, \lambda)_{\hH},
\]
and that $(\hH \times H)_{\varphi} = (\hH \times H)_{\lambda} M$.
Then we have
\[
(G, N, \theta)_{\hH} \geqslant_{\ast} (H, M, \varphi)_{\hH}.
\]
\end{thm}
\begin{proof}
By assumption, we have $G_{\chi^{\hH}} = K H_{\lambda^{\hH}}$, $K \cap H_{\lambda^{\hH}} = Z$, and $(\hH \times H_{\lambda^{\hH}})_{\chi} = (\hH \times H_{\lambda^{\hH}})_{\lambda}$.
Note that $(\hH \times H_{\lambda^{\hH}})_{\lambda} = (\hH \times H)_{\lambda}$. We also have $(\hH \times H_{\lambda^{\hH}})_{\chi} = (\hH \times H)_{\chi}$. To see this, suppose $(\sigma, g) \in (\hH \times H)_{\chi}.$ 
Since  $g \in G_{\chi^{\hH}}$ and $G_{\chi^{\hH}} = K H_{\lambda^{\hH}}$, we can write $g = k h$ for some $k \in K$ and $h \in H_{\lambda^{\hH}}$.
As both $g$ and $h$ lie in $H$, we have $k = g h^{-1} \in H$.
Thus, $k \in K \cap H = Z$, and so $g = k h \in H_{\lambda^{\hH}}$.
This shows that $(\hH \times H)_{\lambda} = (\hH \times H)_{\chi}$, and note that this is a subgroup of $(\hH \times H)_{\theta} = (\hH \times H)_{\varphi}$.

Note that $H_{\lambda}$ is a subgroup of $H_{\varphi}$, and $M \cap H_{\lambda} = M_{\lambda} = Z$.
From $(\hH \times H)_{\varphi} = (\hH \times H)_{\lambda} M$, we deduce that $H_{\varphi} = M H_{\lambda}$.
Similarly, we can prove that $N \cap G_{\chi} = K$ and $G_{\theta} = N G_{\chi}$.
More explicitly, since $N$ and $G_{\chi}$ are subgroups of $G_{\theta}$, and
\[
G_{\theta} =  N H_{\varphi} = N M H_{\lambda} = N H_{\chi} \subseteq N G_{\chi},
\]
we have $G_{\theta} = N G_{\chi}$.

Suppose that $(\Pj, \Pj')$ gives
\[
(G_{\chi^{\hH}}, K, \chi)_{\hH} \geqslant_{\ast} (H_{\lambda^{\hH}}, Z, \lambda)_{\hH}.
\]
Let $\widehat{\Pj} = \Ind^{G_{\theta}}_{G_{\chi}, N}(\Pj)$ and $\widehat{\Pj}' = \Ind^{H_{\varphi}}_{H_{\lambda}, M}(\Pj')$. We will show that $(\widehat{\Pj}, \widehat{\Pj}')$ gives
\[
(G, N, \theta)_{\hH} \geqslant_{\ast} (H, M, \varphi)_{\hH}.
\]
Note that $G_{\theta}/N \cong H_{\varphi}/M$ and $G_{\chi}/K \cong H_{\lambda}/Z$ by assumption, and we have previously shown that $G_{\theta}/N \cong G_{\chi}/K$ and $H_{\varphi}/M \cong H_{\lambda}/Z$. Clearly, these four isomorphisms is natural.
By assumption and the description above the theorem, the factor sets of $\Pj$, $\Pj'$, $\widehat{\Pj}$, and $\widehat{\Pj}'$ coincide under the isomorphism
\[
G_{\chi}/K \cong H_{\lambda}/Z \cong G_{\theta}/N \cong H_{\varphi}/M.
\]

Now suppose that $a = (\sigma, g) \in (\hH \times H)_{\lambda}$. Clearly, $g \in H_{\lambda^{\hH}}$, and thus $g$ normalizes $Z$ and $H_{\lambda}$.
Let $\mu'_a \colon H_{\lambda} \to F^{\times}$ and $\hat{\mu}'_a \colon H_{\varphi} \to F^{\times}$ be the functions determined by $\Pj'^a \sim \mu'_a \Pj'$ and $\widehat{\Pj}'^a \sim \hat{\mu}'_a \widehat{\Pj}'$, respectively.
We want to prove that $\mu'_a$ and $\hat{\mu}'_a$ coincide under the isomorphism $H_{\varphi}/M \cong H_{\lambda}/Z$.
Let $\{m_1 = 1, \dots, m_s\}$ be a complete set of representatives of right $Z$-cosets in $M$, such that for any $x \in H_{\varphi}$,
\[
\widehat{\Pj}'(x) = \diag\big( \Pj'(m_1 x m_{w_x(1)}^{-1}), \dots, \Pj'(m_s x m_{w_x(s)}^{-1}) \big) T_{w_x},
\]
where $w_x \in \bS_s$ is defined by $H_{\lambda} m_i x = H_{\lambda} m_{w_x(i)}$ for $1\leq i \leq s$, and $T_{w_x} := T_{w_x, \lambda(1)}$ is the permutation matrix defined above the theorem.
Note that $\{m_1, \dots, m_s\}$ is also a complete set of representatives of right $H_{\lambda}$-cosets in $H_{\varphi}$.
Since $g$ normalizes both $H_{\varphi}$ and $H_{\lambda}$, let $w \in \bS_s$ be defined by $(H_{\lambda} m_i)^g = H_{\lambda} m_{w(i)}$.
It is straightforward to verify that $w_{g x g^{-1}} = w w_x w^{-1}$ for any $x \in H_{\varphi}$.
Let $z_i \in Z$ be such that $m_i^g = z_i m_{w(i)}$ for each $i$.
For any $x \in H_{\varphi}$, we compute:
\begin{align*}
  \widehat{\Pj}'^a(x)&=\widehat{\Pj}'(gxg^{-1})^{\sigma}=\diag\Big(  \Pj'(m_1 gxg^{-1} m_{w_{gxg^{-1}}(1)}^{-1})^{\sigma},\cdots, \Pj'(m_s gxg^{-1} m_{w_{gxg^{-1}}(s)}^{-1})^{\sigma} \Big)T_{w_{gxg^{-1}}} \\
   &=\diag\Big(  \Pj'^a(m_1^gx (m_{w_{gxg^{-1}}(1)}^g)^{-1}),\cdots, \Pj'^a(m_s^g x(m_{w_{gxg^{-1}}(s)}^g)^{-1}) \Big)T_{w_{gxg^{-1}}}  \\
   &\sim\mu'_a(x)\diag\Big(  \Pj'(m_1^gx (m_{w_{gxg^{-1}}(1)}^g)^{-1}),\cdots, \Pj'(m_s^g x(m_{w_{gxg^{-1}}(s)}^g)^{-1}) \Big)T_{w_{gxg^{-1}}}   \\
   &=\mu'_a(x)\diag\Big(  \Pj'(z_1m_{w(1)}x m_{ww_x(1)}^{-1}z_{w_{gxg^{-1}}(1)}^{-1}),\cdots, \Pj'(z_sm_{w(s)}x m_{ww_x(s)}^{-1}z_{w_{gxg^{-1}}(s)}^{-1}) \Big)T_{w_{gxg^{-1}}}  \\
   &=\mu'_a(x)A\diag\Big(  \Pj'(m_{w(1)}x m_{ww_x(1)}^{-1}),\cdots, \Pj'(m_{w(s)}x m_{ww_x(s)}^{-1}) \Big)T_{w_{gxg^{-1}}}A^{-1}  \\
 &\sim\mu'_a(x)\diag\Big(  \Pj'(m_{w(1)}x m_{ww_x(1)}^{-1}),\cdots, \Pj'(m_{w(s)}x m_{ww_x(s)}^{-1}) \Big)T_{w}T_{w_x}T_{w}^{-1}  \\
 &=\mu'_a(x)T_{w}\diag\Big(  \Pj'(m_{1}x m_{w_x(1)}^{-1}),\cdots, \Pj'(m_{s}x m_{w_x(s)}^{-1}) \Big)T_{w_x}T_{w}^{-1}  \\
 &\sim \mu'_a(x)\diag\Big(  \Pj'(m_{1}x m_{w_x(1)}^{-1}),\cdots, \Pj'(m_{s}x m_{w_x(s)}^{-1}) \Big)T_{w_x}=\mu'_a(x)\widehat{\Pj}'(x), 
\end{align*}
where we define $\mu'_a(x) = \mu'_a(x_1)$ if $x = u x_1$ for some $u \in M$ and $x_1 \in H_{\lambda}$, and in the sixth line, $A = \diag(\Pj'(z_1), \dots, \Pj'(z_s))$.
Thus, we have $\hat{\mu}'_a = \mu'_a$, as desired.

Let $\mu_a \colon G_{\chi} \to F^{\times}$ and $\hat{\mu}_a \colon G_{\theta} \to F^{\times}$ be the functions determined by $\Pj^a \sim \mu_a \Pj$ and $\widehat{\Pj}^a \sim \hat{\mu}_a \widehat{\Pj}$, respectively.
Similarly, we can prove that $\mu_a$ and $\hat{\mu}_a$ coincide under the isomorphism $G_{\theta}/N \cong G_{\chi}/K$.
Since $\mu_a$ and $\mu'_a$ coincide under the isomorphism $G_{\chi}/K \cong H_{\lambda}/Z$ by assumption, it follows that $\hat{\mu}_a$ and $\hat{\mu}'_a$ coincide under the isomorphism $G_{\theta}/N \cong H_{\varphi}/M$.
Given that $(\hH \times H)_{\varphi} = (\hH \times H)_{\lambda} M$ and by Lemma \ref{lem:H-tri.ind}, we conclude that
\[
(G, N, \theta)_{\hH} \geqslant (H, M, \varphi)_{\hH}.
\]

Now further assume that $(\Pj, \Pj')$ gives
\(
(G_{\chi^{\hH}}, K, \chi)_{\hH} \geqslant_c (H_{\lambda^{\hH}}, Z, \lambda)_{\hH}.
\)
If $x \in \C_G(N)$, then it is easy to show that $x \in H_{\lambda}$ and centralizes $K$.
Thus, $\Pj(x)$ and $\Pj'(x)$ are associated with the same scalar.
By the definition of $\widehat{\Pj}$ and $\widehat{\Pj}'$, it is straightforward to verify that $\widehat{\Pj}(x)$ and $\widehat{\Pj}'(x)$ are associated with the same scalar.
This proves that
\[
(G, N, \theta)_{\hH} \geqslant_c (H, M, \varphi)_{\hH}.
\]

Now consider the case $\ast = b$.
It remains to show that for any $x \in G_{\theta}$, 
the scalar matrices $\widehat{\Pj}(\Cl_{\langle N, x \rangle}(x)^+)$ and $\widehat{\Pj}'((\Cl_{\langle N, x \rangle}(x) \cap H_{\varphi})^+)$ are associated with the same scalar.
Fix $x \in G_{\theta}$, and let $J = \langle N, x \rangle$.
By computing the $(1,1)$-entry of the matrix
\[
\sum_{y \in \Cl_J(x)} \widehat{\Pj}(y),
\]
we find that $\widehat{\Pj}(\Cl_J(x)^+)$ and $\Pj((\Cl_J(x) \cap G_{\chi})^+)$ are associated with the same scalar.
Similarly, $\widehat{\Pj}'((\Cl_J(x) \cap H_{\varphi})^+)$ and $\Pj'((\Cl_J(x) \cap H_{\lambda})^+)$ are associated with the same scalar.
Since $\Pj((\Cl_J(x) \cap G_{\chi})^+)$ and $\Pj'((\Cl_J(x) \cap H_{\lambda})^+)$ are associated with the same scalar by assumption, the result follows.
\end{proof}

\section{The Dade--Glauberman--Nagao correspondence}
The main theorem of this section strengthens \cite[Theorem 4.1]{FFZ}. 
For the definition and basic properties of the Dade--Glauberman--Nagao (DGN) correspondence, see \cite[Section 1]{Fu24}.

\begin{hyp}\label{hyp:DGN}
Let $N \zg A$ be finite groups, and let $M/N$ be a normal $p$-subgroup of $A/N$. Suppose that $\theta \in \dz(N)$ is $M$-invariant. Let $D$ be a defect group of the unique block of $M$ covering $\bl(\theta)$. 
Assume that $A = A_{\theta^{\hH}}$.
Let $C = \C_N(D)$, and let $\varphi := \theta^{\star_D} \in \dz(C)$ be the DGN correspondent of $\theta$ with respect to $D$.
\end{hyp}

Notice that $A = N \N_A(D)$ and $N \cap \N_A(D) = C$ under the hypothesis. In fact, we have proved that
\begin{equation}\label{equ:DNG1}
  (A, N, \theta)_{\hH} \geqslant (\N_A(D), C, \varphi)_{\hH}
\end{equation}
in \cite{Fu24} (see Section 4.2 there).
Suppose that \eqref{equ:DNG1} is given by $(\Pj, \Pj')$, and let $\nu$ be the isomorphism of character triples corresponding to $(\Pj, \Pj')$.
We will prove that $\nu$ preserves the Brauer correspondence of blocks.
For the reader's convenience, we recall the construction of $(\Pj, \Pj')$; see \cite{Fu24} for more details.

Let $E = \Fp[\theta] = \Fp[\varphi]$.
Let $\mathfrak{A} = \M_{\theta(1)}(E)$ and $\mathfrak{B} = \M_{\varphi(1)}(E)$.
Let $X \colon N \to \GL_{\theta(1)}(E) = \mathfrak{A}^{\times}$ and $X' \colon C \to \GL_{\varphi(1)}(E) = \mathfrak{B}^{\times}$ be group representations of $N$ and $C$ affording $\theta$ and $\varphi$, respectively, both realized over $E$.
In fact, $\mathfrak{A}$ is a Dade $D$-algebra under the action spanned by $X(n)^d = X(n^d)$ for $n \in N$ and $d \in D$, and $\mathfrak{B}$ is the $D$-Brauer quotient of $\mathfrak{A}$ with $X'(c) = \Br_D(X(c))$ for any $c \in C$.
There is a unique group homomorphism $\rho \colon D \to \mathfrak{A}^{\times}$ realizing the $D$-action by conjugation, and we let $D_1 = \rho(D)$.
There exists a group homomorphism $\phi \colon \N_{\mathfrak{A}^{\times}}(D_1) \to \mathfrak{B}^{\times}$ (not unique in general) such that $\phi(x) = \Br_D(x)$ for $x \in (\mathfrak{A}^D)^{\times}$.
Choose $\Pj \colon A_{\theta} \to \GL_{\theta(1)}(E)$ to be a projective representation associated with $\theta$ and realized over $E$, and define
\[
\Pj' \colon \N_A(D)_{\varphi} \to \GL_{\varphi(1)}(E), \quad g \mapsto \phi(\Pj(g)).
\]
Then the pair $(\Pj, \Pj')$ gives \eqref{equ:DNG1}.

Let $\hat{\theta} \in \IBr(M)$ and $\hat{\varphi} \in \IBr(\N_M(D))$ be the unique characters lying over $\theta$ and $\varphi$, respectively.
Recall that $\hat\theta_N = \theta$ and $\hat\varphi_C = \varphi$.

\begin{lem}\label{lem:DGN1}
Keeping the notation above, we may assume that $\Pj$ and $\Pj'$ are associated with $\hat\theta$ and $\hat\varphi$, respectively.  
\end{lem}
\begin{proof}
Since $\hat\theta_N = \theta$ and $\hat\theta$ is the unique character lying over $\theta$, by the fundamental theorem of Galois theory we have $\Fp[\theta] = \Fp[\hat\theta]$.
Since any projective representation of $A_{\theta}$ associated with $\hat\theta$ is also associated with $\theta$, by adjusting scalars in $E^{\times}$, we may assume that $\Pj$ is associated with $\hat\theta$. This means that $\Pj_M$ is a group representation and
\[
\Pj(g) \Pj(h) = \Pj(gh), \quad \Pj(h) \Pj(g) = \Pj(hg)
\]
for all $g \in A_{\theta}$ and $h \in M$.
Since $\Pj'(g) = \phi(\Pj(g))$ for all $g \in \N_A(D)_{\varphi}$ and $\phi$ is a group homomorphism, it follows that $\Pj'_{\N_M(D)}$ is a group homomorphism and
\[
\Pj'(g) \Pj'(h) = \Pj'(gh), \quad \Pj'(h) \Pj'(g) = \Pj'(hg)
\]
for all $g \in \N_A(D)_{\varphi}$ and $h \in \N_M(D)$. This implies that $\Pj'$ is associated with $\hat\varphi$, and completes the proof.
\end{proof}

The following theorem is of vital importance for studying the block isomorphism of $\hH$-triples arising from the DGN correspondence.
It strengthens the main result of \cite{Fu24}, and the proof we present adopts a group algebra perspective.

\begin{thm}\label{thm:DGN1}
Keeping the notation above, the isomorphism $\nu$ of character triples corresponding to $(\Pj, \Pj')$ satisfies that for any $G \in \mS(A_{\theta}, M)$ and $\chi \in \IBr(G \which \theta)$, we have
\[
\bl(\nu_G(\chi))^G = \bl(\chi).
\]
\end{thm}
\begin{proof}
Let $G$ and $\chi$ be as in the theorem, and fix them.
Consider the group algebra $FG$, and let the Brauer homomorphism
\[
\Br_D \colon (FG)^D \to F\N_G(D)
\]
be defined by
\[
\Br_D(\Cl_D(x)^+) = 
\begin{cases}
x, & x \in \C_G(D), \\
0, & x \in G \setminus \C_G(D),
\end{cases}
\]
where $(FG)^D$ is the subalgebra of $D$-invariant elements in $FG$.
Let $e_{\theta} \in \Z(FN)$ and $e_{\varphi} \in \Z(FC)$ be the block idempotents associated with $\bl(\theta)$ and $\bl(\varphi)$, respectively.
Note that $e_{\theta} \in \Z(FG)$ since $\theta$ is $G$-invariant.
Also, $e_{\varphi} \in \Z(F\N_G(D))$.
Since $e_{\varphi} = \Br_D(e_{\theta})$ by the definition of the DGN correspondence, the Brauer homomorphism restricts to
\[
\Br_D \colon (FG e_{\theta})^D \to F\N_G(D) e_{\varphi},
\]
and we still denote it by the same symbol.

Let $f \in \Z(FG)$ be the block idempotent of $\bl(\chi)$.
Note that $f$ belongs to $e_{\theta}$, meaning that $e_{\theta} f = f$.
We will prove that $\Br_D(f)$ is a block idempotent of $F\N_G(D)$ associated with the block $\bl(\nu_G(\chi))$.
Then by \cite[Chapter, Theorem 3.5]{NagTsu}, we have
\[
\bl(\nu_G(\chi))^G = \bl(\chi),
\]
which proves the theorem.

Note that $FN e_{\theta}$ is a full matrix algebra over the field $F$, and 
\[
FG e_{\theta} = FN e_{\theta} \otimes \C_{FG e_{\theta}}(FN e_{\theta})
\]
(see \cite[Chapter V, Theorem 7.2]{NagTsu}).
Now $EN e_{\theta}$ is a subalgebra (over $E$) of $FN e_{\theta}$, and there is an isomorphism of algebras from $EN e_{\theta}$ to $\mathfrak{A}$ given by $n e_{\theta} \mapsto X(n)$ for $n \in N$.
For any $g \in G$, let $s_g \in (EN e_{\theta})^{\times}$  be the element corresponding to $\Pj(g)$ under this isomorphism.
Since $FN e_{\theta}$ is $F$-spanned by $\{n e_{\theta} \which n \in N\}$, the conjugation actions of $s_g$ and $g$ on $FN e_{\theta}$ agree.

Let $\bar G = G/N = \N_G(D)/C$.
Note that $s_g^{-1} g \in \C_{FG e_{\theta}}(FN e_{\theta})$ is independent of the choice of $g \in G$ in its $N$-coset. For any $ g \in  G$, let $u_{\bar g} = s_g^{-1} g$.
It is straightforward to verify that $\C_{FG e_{\theta}}(FN e_{\theta})$ is a generalized group ring with $F$-basis $\{u_{\bar g} \which \bar g \in \bar G\}$ and multiplication given by
\[
u_{\bar g} u_{\bar h} = \bar\alpha(\bar g, \bar h)^{-1} u_{\bar g \bar h},
\]
where $\bar\alpha$ is the factor set associated with $\Pj$ (see also the proof of \cite[Chapter V, Theorem 7.2]{NagTsu}).
Similarly, we have
\[
F\N_G(D) e_{\varphi} = FC e_{\varphi} \otimes \C_{F\N_G(D) e_{\varphi}}(FC e_{\varphi}).
\]
Since $ECe_{\varphi}$ is isomorphic to $\mathfrak{B}$, for $g \in \N_G(D)$, let $s'_g\in (ECe_{\varphi})^{\times}$ be the element corresponding to $\Pj'(g)$, and let $v_{\bar g} = (s'_g)^{-1} g$.
Then $\C_{F\N_G(D) e_{\varphi}}(FC e_{\varphi})$ is a generalized group ring with $F$-basis $\{v_{\bar g} \which \bar g \in \bar G\}$ and is isomorphic to $\C_{FG e_{\theta}}(FN e_{\theta})$ via the map sending $v_{\bar g}$ to $u_{\bar g}$.

It is straightforward to verify that if $g \in \C_G(D)$, then $s_g$ is $D$-invariant and $s'_g = \Br_D(s_g)$.
Note that $\C_G(D) \in \mS(\N_G(D), C)$.
Let $\mathcal{T}$ be a complete set of representatives of $C$-cosets in $\N_G(D)$ with $1 \in \mathcal{T}$.
Then $\mathcal{T}$ is also a complete set of representatives of $N$-cosets in $G$.

Since $f$ is a primitive central idempotent of $FG e_{\theta}$ and\[
   \Z(FGe_{\theta}) =\Z(FNe_{\theta})\otimes\Z(\C_{FGe_{\theta}}(FNe_{\theta})) =1\otimes \Z(\C_{FGe_{\theta}}(FNe_{\theta})),
\]
it follows that $f$ is a primitive central idempotent of $\C_{FG e_{\theta}}(FN e_{\theta})$. Write
\[
f = \sum_{g \in \mathcal{T}} \lambda_g u_{\bar g}, \quad \lambda_g \in F.
\]
Let
\[
f' = \sum_{g \in \mathcal{T}} \lambda_g v_{\bar g} \in \C_{F\N_G(D) e_{\varphi}}(FC e_{\varphi}).
\]
Since $\C_{F\N_G(D) e_{\varphi}}(FC e_{\varphi})$ is isomorphic to $\C_{FG e_{\theta}}(FN e_{\theta})$, we conclude that $f'$ is a primitive central idempotent of $\C_{F\N_G(D) e_{\varphi}}(FC e_{\varphi})$, and hence also a primitive central idempotent of $F\N_G(D) e_{\varphi}$ by analogous arguments.
In fact, $f'$ is the block idempotent of $\bl(\nu_G(\chi))$, as can be seen directly from the definition.

It remains to prove that $f' = \Br_D(f)$.
Since $f' = \sum_{g \in \mathcal{T}} \lambda_g v_{\bar g}$ is a block idempotent of $\N_G(D)$, by \cite[Chapter V, Theorem 2.8]{NagTsu}, we have $\lambda_g = 0$ if $g \notin \C_G(D)$.
Let $\mathcal{T}_1 = \mathcal{T} \cap \C_G(D)$. Then
\[
f = \sum_{g \in \mathcal{T}_1} \lambda_g u_{\bar g}.
\]
Thus,
 \begin{align*}
   \Br_D(f)&=\Br_D\Big(\sum_{g\in\mathcal{T}_1}\lambda_gu_{\bar g}\Big) =\sum_{g\in\mathcal{T}_1}\lambda_g\Br_D (s_g^{-1}g) =\sum_{g\in\mathcal{T}_1}\lambda_g\Br_D (s_g)^{-1}g  \\
   &
   =\sum_{g\in\mathcal{T}_1}\lambda_g(s'_g)^{-1}g=f'.
 \end{align*}
This completes the proof of the theorem.
\end{proof}

\begin{thm}\label{thm:DGN-cor}
Assume Hypothesis \ref{hyp:DGN}.
For any $G \in \mS(A_{\theta}, M),$  
then there exists an $\big(\hH \times \N_A(D)\big)_{\theta}$-equivariant bijection
\begin{equation}\label{equ:DGN-cor2}
  \nu_G \colon \IBr(G \which \theta, |D|) \to \IBr(\N_G(D) \which \varphi, |D|)
\end{equation}
such that 
\[
(A_{\chi^{\hH}}, G, \chi)_{\hH} \geqslant_b (\N_A(D)_{\chi^{\hH}}, \N_G(D), \Delta(\chi))_{\hH}
\]
for any $\chi \in \IBr(G \which \theta, |D|)$.
\end{thm}
\begin{proof}
In \cite{Fu24}, we proved that $A = N \N_A(D)$ and $C = N \cap \N_A(D)$.
Let $\hat{\theta} \in \IBr(M)$ and $\hat{\varphi} \in \IBr(\N_M(D))$ be the unique characters lying over $\theta$ and $\varphi$, respectively.

Let $(\Pj, \Pj')$ be as described above.
We have shown that $(\Pj, \Pj')$ gives
\[
(A, M, \hat\theta)_{\hH} \geqslant (\N_A(D), \N_M(D), \hat\varphi)_{\hH},
\]
and we let $\nu$ be the isomorphism of character triples corresponding to $(\Pj, \Pj')$.
Since $\nu$ preserves vertices of characters by \cite[Theorem A]{Fu24}, and by Lemma \ref{lem:H-tri.B}, it suffices to prove that
\begin{equation}\label{equ:DGN2}
  (A, M, \hat\theta)_{\hH} \geqslant_b (\N_A(D), \N_M(D), \hat\varphi)_{\hH}.
\end{equation}

Let $\omega_c$ be the scalar associated with $\Pj(c)$ for $c \in \C_A(M) \subseteq \C_{A_{\theta}}(D)$.
Since $\Pj'(c) = \Br_D(\Pj(c))$, it follows that $\Pj(c)$ and $\Pj'(c)$ are associated with the same scalar $\omega_c$.
This proves that
\[
(A, M, \hat\theta)_{\hH} \geqslant_c (\N_A(D), \N_M(D), \hat\varphi)_{\hH}.
\]
Then \eqref{equ:DGN2} follows from Theorem \ref{thm:2.two.equiv.}(2) and Theorem \ref{thm:DGN1}.
\end{proof}

\section{The inductive NAW (BNAW)  condition}
In this section, we introduce the inductive NAW (resp. BNAW) conditions for finite non-abelian simple groups. We prove that Conjecture \ref{ConjA} (resp. Conjecture \ref{ConjB}) holds for central extensions of a direct product of isomorphic non-abelian simple groups satisfying the inductive NAW (resp. BNAW) condition (see Theorem \ref{thm:ibGAW.con.}).

Recall that the \emph{universal $p'$-covering group} of a perfect group $L$ is the maximal perfect central extension of $L$ by an abelian $p'$-group (see \cite[Appendix B]{Navarro:McKay} for universal covering groups).

\begin{defi}\label{def:ibGAW}
Let $L$ be a finite non-abelian simple group of order divisible by $p$, and let $S$ be the universal $p'$-covering group of $L$.
We say that the inductive NAW condition (resp. inductive BNAW condition) holds for $L$ at $p$ if Conjecture \ref{ConjA} (resp. Conjecture \ref{ConjB}) holds for $S$ at $p$.
\end{defi}

\begin{lem}\label{lem:ibGAW}
Let $m \geq 1$ be an integer, and let $L$ be a finite non-abelian simple group of order divisible by $p$ satisfying the inductive NAW (resp. BNAW) condition.
Let $S$ be the universal $p'$-covering group of $L$.
Then Conjecture \ref{ConjA} (resp. Conjecture \ref{ConjB}) holds for $S^m$.
\end{lem}
\begin{proof}
Since every $\widetilde Q \in \Rad(S^m)$ can be written as $\widetilde Q = Q_1 \times \cdots \times Q_m$ with $Q_i \in \Rad(S)$ for $1\leq i \leq m$, it follows from \cite[Lemma 2.3 (b)]{NT11} that
\begin{align*}
  \IBr(S^m) &= \IBr(S) \times \cdots \times \IBr(S) \quad \text{($m$ times)}, \quad \text{and} \\
  \W^{\circ}(S^m)/\sim_{S^m} &= \W^{\circ}(S)/\sim_S \times \cdots \times \W^{\circ}(S)/\sim_S \quad \text{($m$ times)}.
\end{align*}
By assumption, there exists an $\hH \times \Aut(S)$-equivariant bijection 
\[
\Omega_S \colon \IBr(S) \to \W^{\circ}(S)/\sim_S
\]
such that $(A_{\theta^{\hH}},S, \theta)_{\hH} \geqslant_{\ast} (\N_A(Q)_{\varphi^{\hH}},\N_S(Q), \varphi)_{\hH},$ where $A=S\semi\Aut(S), \theta \in \IBr(S)$ and $(Q, \varphi) \in \Omega_S(\theta)$, where $\ast \in \{c, b\}$ depending on our assumption.

Define the map
\begin{align*}
  \Omega_{S^m} \colon \IBr(S^m) &\to \W(S^m)/\sim_{S^m}, \\
  \theta_1 \times \cdots \times \theta_m &\mapsto \Omega_S(\theta_1) \times \cdots \times \Omega_S(\theta_m),
\end{align*}
where $\theta_1, \dots, \theta_m \in \IBr(S)$.
By \cite[Lemma 10.24]{Navarro:McKay}, we have $\Aut(S^m) = \Aut(S) \wr \bS_m$, so it is straightforward to verify that $\Omega_{S^m}$ is an $\hH \times \Aut(S^m)$-equivariant bijection.

It is convenient to identify $S^m\semi \Aut(S^m)$ with $A \wr \bS_m.$
Following the last paragraph of the proof of \cite[Lemma 7.5]{FFZ}, and uses Lemmas \ref{lem:H-tri1}--\ref{lem:H-tri3} in place of the lemmas used there,
we have \[((A\wr \bS_m)_{\tilde\theta^{\hH}},S^m, \tilde\theta)_{\hH} \geqslant_{\ast} (\N_{A\wr \bS_m}(\widetilde{Q})_{\tilde{\varphi}^{\hH}},\N_{S^m}(\widetilde Q), \tilde{\varphi})_{\hH}\] for any $\tilde\theta \in \IBr(S^m)$ and $(\widetilde Q, \tilde\varphi) \in \Omega_{S^m}(\tilde\theta)$. 
The Lemma then follows by Theorem \ref{thm:butterfly} and Proposition \ref{prop:lieover2}.
\end{proof}

\begin{prop}\label{prop:ibGAW1}
Let $K$ be a finite perfect group such that $Z := \Z(K)$ is a $p'$-group and $K/Z$ is a direct product of isomorphic non-abelian simple groups satisfying the inductive NAW (resp. BNAW) condition.
Then Conjecture \ref{ConjA} (resp. Conjecture \ref{ConjB}) holds for $K$.
\end{prop}
\begin{proof}
Let $K/Z \cong L^m$, and let $S$ be the universal $p'$-covering group of $L$.
Let $\epsilon \colon S^m \to K$ be the universal $p'$-central extension of $K$.
Since every automorphism $\phi$ of $K$ lifts to a unique automorphism $\hat\phi$ of $S^m$, we can regard $\Aut(K)$ as a subgroup of $\Aut(S^m)$.
In fact, $\Aut(K) = \Aut(S^m)_{\ker(\epsilon)}$.
Let $1_{\ker(\epsilon)}$ be the trivial character of $\ker(\epsilon)$.
Then we can identify $\IBr(K)$ with the subset $\IBr(S^m \which 1_{\ker(\epsilon)})$ of $\IBr(S^m)$.

Since there is a natural correspondence between $\Rad(K)$ and $\Rad(S^m)$ by taking the normal Sylow $p$-subgroup of $\epsilon^{-1}(Q)$ for $Q \in \Rad(K)$, we can regard $\W^{\circ}(K)$ as a subset of $\W^{\circ}(S^m)$.
In fact, $\W^{\circ}(K) = \W^{\circ}(S^m \which 1_{\ker(\epsilon)})$, where $\W^{\circ}(S^m \which 1_{\ker(\epsilon)})$ consists of weights $(Q, \delta)$ of $S^m$ such that $\delta$ lies over $1_{\ker(\epsilon)}$.

By assumption, there exists an $\hH \times \Aut(S^m)$-equivariant bijection
\[
\Omega_{S^m} \colon \IBr(S^m) \to \W(S^m)/\sim_{S^m}
\]
such that for any $\theta \in \IBr(S^m\which 1_{\ker(\epsilon)})$ and $(Q, \varphi) \in \Omega_{S^m}(\theta)$, we have
\begin{equation}\label{equ:iBGAW1}
  (S^m\semi\Aut(S^m)_{\ker(\epsilon),\theta^{\hH}},S^m, \theta)_{\hH} \geqslant_{\ast} (\N_{S^m\semi\Aut(S^m)_{\ker(\epsilon)}}(Q)_{\varphi^{\hH}},\N_{S^m}(Q), \varphi)_{\hH}.
\end{equation}
The bijection $\Omega_{S^m}$ restricts to an $\hH \times \Aut(K)$-equivariant bijection
\[
\Omega_K \colon \IBr(K) \to \W(K)/\sim_K.
\]

Write \[\overline{S^m\semi \Aut(S^m)_{\ker\epsilon,\theta^{\hH}}}=S^m\semi \Aut(S^m)_{\ker\epsilon,\theta^{\hH}}/\ker\epsilon\]
and identify it with $K\semi\Aut(K)_{\bar\theta^{\hH}}.$
Since every non-trivial automorphism $\hat\phi \in \Aut(S^m)_{\ker(\epsilon)}$ descends to a non-trivial automorphism $\phi \in \Aut(K)$, Lemma \ref{lem:H-tri.desc.} applies to (\ref{equ:iBGAW1}) and yields
\[
(K\semi\Aut(K)_{\bar\theta^{\hH}}, K, \bar\theta)_{\hH} \geqslant_{\ast} (\N_{K\semi\Aut(K)}(\bar Q)_{\bar\varphi^{\hH}}, \N_K(\bar Q), \bar\varphi)_{\hH}.
\]
The Proposition then follows by Theorem \ref{thm:butterfly} and Proposition \ref{prop:lieover2}.
\end{proof}

\begin{thm}\label{thm:ibGAW.con.}
Let $Z$ be a cyclic central $p'$-subgroup of a finite group $K$.
Suppose that $K/Z$ is either a direct product of isomorphic non-abelian simple groups satisfying the inductive NAW (resp. BNAW) condition, or a $p'$-group.
Then Conjecture \ref{ConjA} (resp. Conjecture \ref{ConjB}) holds for the group $K$.
\end{thm}
\begin{proof}
The theorem holds trivially if $K/Z$ is a $p'$-group.
Now suppose that $K/Z$ is a direct product of isomorphic non-abelian simple groups satisfying the inductive NAW (resp. BNAW) condition.
Let $K_1$ be the commutator subgroup of $K$ and $Z_1 = Z \cap K_1 = \Z(K_1)$. Note that $K_1$ is perfect and $K = K_1 Z$.
Since $K_1/Z_1 \cong K/Z$, by Proposition \ref{prop:ibGAW1}, Conjecture \ref{ConjA} (resp. Conjecture \ref{ConjB}) holds for $K_1$.

Thus, for any $K \zg A$, there exists an $\hH \times A$-equivariant bijection $\Omega_{K_1} \colon \IBr(K_1) \to \W(K_1)/\sim_{K_1}$ such that for any $\theta \in \IBr(K_1)$ and $(Q, \varphi) \in \Omega_{K_1}(\eta)$,
\begin{equation}\label{equ:ibGAW1}
  (A_{\theta^{\hH}}, K_1, \theta)_{\hH} \geqslant_{\ast} (\N_A(Q)_{\varphi^{\hH}}, \N_{K_1}(Q), \varphi)_{\hH}.
\end{equation}

By \cite[Lemma 2.2]{NT11}, every character in $\IBr(K)$ can be written uniquely as $\theta \cdot \lambda$ for some $\theta \in \IBr(K_1)$ and $\lambda \in \IBr(Z)$ such that $\theta$ and $\lambda$ lie over the same irreducible Brauer character of $Z_1$.
We define $\Omega_K(\theta \cdot \lambda)$ to be the $K$-conjugacy class containing the $p$-weight $(Q, \varphi \cdot \lambda)$ of $K$, where $(Q, \varphi) \in \Omega_{K_1}(\theta)$. (Note that $\Rad(K) = \Rad(K_1)$.)
It is straightforward to verify that the map $\Omega_K \colon \IBr(K) \to \W(K)/\sim_K$ is a well-defined $\hH \times A$-equivariant bijection.

It remains to prove that
\begin{equation}\label{equ:ibGAW2}
  \big(A_{(\theta \cdot \lambda)^{\hH}}, K, \theta \cdot \lambda\big)_{\hH} \geqslant_{\ast} 
  \big(\N_A(Q)_{(\varphi \cdot \lambda)^{\hH}}, \N_K(Q), \varphi \cdot \lambda\big)_{\hH}.
\end{equation}
Suppose that \eqref{equ:ibGAW1} is given by $(\Pj,\Pj')$ and $\nu$ is the isomorphism of character triples corresponding to $(\Pj,\Pj')$.
Since $\theta \cdot \lambda$ (resp. $\varphi \cdot \lambda$) is the unique character in $\IBr(K \which \theta)$ (resp. $\IBr(\N_K(Q) \which \varphi)$) lying over $\lambda$, by Theorem \ref{thm:2.two.equiv.} (1), we have $\nu_K(\theta \cdot \lambda) = \varphi \cdot \lambda$.
Then  \eqref{equ:ibGAW2} follows from Lemma \ref{lem:H-tri.B}.
This completes the proof of the theorem.
\end{proof}

\section{The reduction}
In this section, we present our final reduction, which leads directly to the proof of Theorem C by taking $Z = 1$.

\begin{thm}[reduction]
Let $Z \zg G$ be finite groups, and let $\lambda \in \dz(Z)$ be $G$-invariant. Assume that $G$ is normally embedded in a finite group $A$ with $Z$ normal in $A$.
Let
\[
\rw^{\circ}(G\which \lambda) = \left\{ (S, \varrho) \Bwhich S/Z \in \Rad(G/Z),\ \varrho \in \IBr(\N_G(S) \which \lambda, |S/Z|) \right\}.
\]
Assume that the inductive NAW (resp. BNAW) condition holds for every non-abelian simple group involved in $G/Z$ that has order divisible by $p$.

Then there exists an $(\hH \times A)_{\lambda}$-equivariant bijection
\(
f \colon  \IBr(G \which \lambda) \to \rw^{\circ}(G\which \lambda)/\sim_G
\)
such that for any $\psi\in\IBr(G\which\lambda)$ and $(S, \varrho) \in f(\psi)$, we have
\[
(A_{\psi^{\hH}}, G, \psi)_{\hH} \geqslant_{\ast} (\N_A(S)_{\delta^{\hH}}, \N_G(S), \varrho)_{\hH},
\]
where  $\ast \in \{c,b\}$ depends on our assumption.
\end{thm}
\begin{proof}
The proof of this theorem follows the proof of \cite[Theorem 8.2]{FFZ} with some improvements.
All notations are kept, except that the group $\mathscr{A}$ is replaced by $(\hH\times A)_{\lambda}$ here.

We proceed by induction on $|G/Z|$.
If $|G/Z| = 1$, the theorem holds trivially, so we assume $|G/Z| > 1$.
Without loss of generality, we may assume $A = A_{\lambda^{\hH}}$.
For any $(S, \varrho) \in \rw^{\circ}(G, \lambda)$, any vertex $V$ of $\delta$ must intersect $Z$ trivially. This follows from Fong's theorem \cite[Chapter V, Theorem 5.16(ii)]{NagTsu} and \cite[Chapter V, Theorem 1.9(i)]{NagTsu}.

Applying Theorem \ref{thm:centrali.} to the $\hH$-triple $(A,Z,\lambda)_{\hH}$, we may assume that $Z$ is a central $p'$-subgroup of $G$ and $\lambda$ is a faithful linear character of $Z$.
Since there is a natural bijection $\Rad(G) \to \Rad(G/Z)$, $S \mapsto SZ/Z$, the set $\rw^{\circ}(G\which \lambda)$ corresponds naturally to the set
\[
\W^{\circ}(G\which \lambda) = \left\{ (S, \varrho) \bwhich S \in \Rad(G),\ \varrho \in \IBr(\N_G(S) \which \lambda, |S|) \right\}.
\]

In the proof of \cite[Theorem 8.2]{FFZ}, we constructed an $(\hH \times A)_{\lambda}$-equivariant bijection
\[
f \colon  \IBr(G \which \lambda) \to  \W^{\circ}(G\which \lambda)/\sim_G.
\]
Let $\psi\in\IBr(G\which\lambda)$ and $(S,\varrho)\in f(\psi).$  
As in the proof of Proposition \ref{prop:lieover2}, we have $\N_A(S)_{\psi^{\hH}} = \N_A(S)_{\varrho^{\hH}}$, 
$G \cap \N_A(S)_{\varrho^{\hH}} = \N_G(S)$,  $A_{\psi^{\hH}} = G \N_A(S)_{\varrho^{\hH}}$, and $\big(\hH \times \N_A(S)\big)_{\psi} = \big(\hH \times \N_A(S)\big)_{\varrho}$.

It remains to show that
\[
(A_{\psi^{\hH}}, G, \psi)_{\hH} \geqslant_{\ast} (\N_A(S)_{\varrho^{\hH}}, \N_G(S), \varrho)_{\hH}.
\]
To prove this, we trace the construction of the map $f$.
Let $K/Z$ be a minimal normal subgroup of $A/Z$ contained in $G/Z$.

If $K/Z$ is a $p$-group, then $K = K_p \times Z$, where $K_p$ is the normal Sylow $p$-subgroup of $K$.
Let $\bar{A} = A/K_p$. The map $f$ is obtained by applying induction to $(\bar{A}, \bar{G}, \bar{Z}, \bar{\lambda})$; see the proof of \cite[Theorem 8.2]{FFZ} for details.
By induction,
\[
(\bar{A}_{\bar{\psi}^{\hH}}, \bar{G}, \bar{\psi})_{\hH} \geqslant_{\ast} (\N_{\bar{A}}(\bar{S})_{\bar{\varrho}^{\hH}}, \N_{\bar{G}}(\bar{S}), \bar{\varrho})_{\hH},
\]
where $\ast \in \{c, b\}$ depends on our assumption.
Since $\overline{A_{\psi^{\hH}}} = \bar{A}_{\bar{\psi}^{\hH}}$ and $\overline{\N_A(S)_{\varrho^{\hH}}} = \N_{\bar{A}}(\bar{S})_{\bar{\varrho}^{\hH}}$, Lemma \ref{lem:H-tri.lift} gives
\[
(A_{\psi^{\hH}}, G, \psi)_{\hH} \geqslant_{\ast} (\N_A(S)_{\varrho^{\hH}}, \N_G(S), \varrho)_{\hH}.
\]

Now assume $K/Z$ is not a $p$-group, so $\bO_p(G) = 1$.
Since $K/Z$ is characteristically simple, it is a direct product of isomorphic simple groups.
Thus, $K/Z$ is either a direct product of isomorphic non-abelian simple groups satisfying the inductive NAW (resp. BNAW) condition, or a $p'$-group.
We now recall the steps linking $\psi$ and $\delta$ from the proof of \cite[Theorem 8.2]{FFZ}.

Let $\eta \in \IBr(K)$ be an irreducible constituent of $\psi_K$, and let $\varphi \in \IBr(G_{\eta} \which \eta)$ be the Clifford correspondent of $\psi$.

By Theorem \ref{thm:ibGAW.con.}, there exists an $\hH \times A$-equivariant bijection $\Omega \colon \IBr(K) \to \W(K)/\sim_K$ satisfying Conjecture \ref{ConjA} (resp. Conjecture \ref{ConjB}).
Let $\eta\in\IBr(K)$ and $(Q, \delta) \in \Omega(\eta)$. Then
\[
(A_{\eta^{\hH}}, K, \eta)_{\hH} \geqslant_{\ast} (\N_A(Q)_{\delta^{\hH}}, \N_K(Q), \delta)_{\hH},
\]
and we denote the corresponding isomorphism by $\nu^{(\eta, Q)}$.
Note that $\N_A(Q)_{\eta^{\hH}}=\N_A(Q)_{\delta^{\hH}}.$
Let $U_{\eta, Q} = \N_{G_{\eta}}(Q)$. 
Now $\nu^{(\eta, Q)}_{G_{\eta}} \colon \IBr(G_{\eta} \which \eta) \to \IBr(U_{\eta,Q} \which \delta)$ is a bijection,  and let $\varphi' = \nu^{(\eta, Q)}_{G_{\eta}}(\varphi)$.
By Lemma \ref{lem:H-tri.B}, we have
\begin{equation}\label{equ:red4}
  (A_{\varphi^{\hH}}, G_{\eta}, \varphi)_{\hH} \geqslant_{\ast} (\N_A(Q)_{\varphi'^{\hH}}, U_{\eta,Q}, \varphi')_{\hH},
\end{equation}
noting that $A_{\varphi^{\hH}} \subseteq A_{\eta^{\hH}}$, $\N_A(Q)_{\varphi'^{\hH}} \subseteq \N_A(Q)_{\delta^{\hH}}$ and $\N_A(Q)_{\varphi'^{\hH}}=\N_A(Q)_{\varphi^{\hH}}.$

Let $\overline{\N_A(Q)} = \N_A(Q)/Q$.
Then $\overline{\delta}$ is a defect-zero character of $\overline{\N_K(Q)}$ that is $\overline{U_{\eta,Q}}$-invariant, and $U_{\eta,Q}$ is normal in $\N_A(Q)_{\delta^{\hH}}$ by direct computation.
Applying induction to \[\big(\overline{\N_A(Q)_{\delta^{\hH}}}, \overline{U_{\eta,Q}}, \overline{\N_K(Q)}, \overline{\delta} \big),\]
we have an $(\hH\times A)_{\eta,Q}$-equivariant bijection
\[\bar{f}_3^{(\eta,Q)}:\IBr(\overline{U_{\eta,Q}}\which\bar\delta)\ra \rw^{\circ}(\overline{U_{\eta,Q}}\which \bar\delta)/\sim_{\overline{U_{\eta,Q}}},\]
such that for the $\overline{\varphi'}\in \IBr(\overline{U_{\eta,Q}}\which\bar\delta)$ and $(\overline{E},\bar\gamma)\in \bar{f}_3^{(\eta,Q)}(\overline{\varphi'}),$ we have
\begin{equation}\label{equ:red1}
  (\overline{\N_A(Q)}_{{\overline{\varphi'}}^{\hH}}, \overline{U_{\eta,Q}}, \overline{\varphi'})_{\hH} \geqslant_{\ast}  
  (\overline{\N_A(Q, E)}_{\bar\gamma^{\hH}}, \N_{\overline{U_{\eta,Q}}}(\overline{E}), \bar\gamma)_{\hH},
\end{equation}
where $\overline{E}/\overline{\N_K(Q)} \in \Rad(\overline{U_{\eta,Q}}/\overline{\N_{K}(Q)})$ and $\bar\gamma \in \IBr(\N_{\overline{U_{\eta,Q}}}(\overline{E}) \which \overline{\delta}, |\overline{E}/\overline{\N_K(Q)}|)$.
Note that $\overline{\N_A(Q)_{\delta^{\hH}}}
=\overline{\N_A(Q)}_{\overline{\delta}^{\hH}}, \overline{\N_A(Q)}_{{\overline{\varphi'}}^{\hH}}\subseteq \overline{\N_A(Q)}_{\overline{\delta}^{\hH}}$ and $\overline{\N_A(Q, E)}_{\bar\gamma^{\hH}}\subseteq \overline{\N_A(Q,E)}_{\overline{\delta}^{\hH}}.$


Let $\overline{F}$ be a defect group of the unique block of $\overline{E}$ covering $\bl(\overline{\delta})$, and let $\overline{C_F} = \C_{\overline{\N_K(Q)}}(\overline{F})$.
Let $(\overline{\delta})^{\star_{\overline{F}}} \in \IBr(\overline{C_F})$ be the DGN correspondent of $\overline{\delta}$ with respect to $\overline{F}$.
Since $\overline{E}$ is normal in $\overline{\N_A(Q, E)}_{\overline{\delta}^{\hH}}$ and $\overline{\N_A(Q, E, F)}_{\overline{\delta}^{\hH}} = \overline{\N_A(F)}_{\overline{\delta}^{\hH}}$ (as $E = F \N_K(Q)$ and $Q = K \cap F$), Theorem \ref{thm:DGN-cor} gives a bijection
\[
\bar f_4^{(\eta, Q, E, F)} \colon \IBr(\N_{\overline{U_{\eta,Q}}}(\overline{E}) \which \overline{\delta}, |\overline{F}|) \to \IBr(\N_{\overline{U_{\eta,Q}}}(\overline{F}) \which (\overline{\delta})^{\star_{\overline{F}}}, |\overline{F}|),
\]
such that
\begin{equation}\label{equ:red2}
  \big(\overline{\N_A(Q, E)}_{\bar\gamma^{\hH}}, \N_{\overline{U_{\eta,Q}}}(\overline{E}), \bar\gamma\big)_{\hH} \geqslant_{\ast} \big(\overline{\N_A(F)}_{\bar\zeta^{\hH}}, \N_{\overline{U_{\eta,Q}}}(\overline{F}), \bar\zeta\big)_{\hH},
\end{equation}
where $\bar\zeta = \bar f_4^{(\eta, Q, E, F)}(\bar\gamma)$. Note that $\overline{\N_A(Q, E)}_{\bar\gamma^{\hH}}\subseteq \overline{\N_A(Q, E)}_{\overline{\delta}^{\hH}}$ and $\overline{\N_A(F)}_{\bar\zeta^{\hH}}\subseteq \overline{\N_A(F)}_{\overline{\delta}^{\hH}}.$

Applying Lemma \ref{lem:H-tri.tran.} to \eqref{equ:red1} and \eqref{equ:red2}, we get
\begin{equation}\label{equ:red3}
  (\overline{\N_A(Q)}_{{\overline{\varphi'}}^{\hH}}, \overline{U_{\eta,Q}}, \overline{\varphi'})_{\hH} \geqslant_{\ast} \big(\overline{\N_A(F)}_{\bar\zeta^{\hH}}, \N_{\overline{U_{\eta,Q}}}(\overline{F}), \bar\zeta\big)_{\hH}.
\end{equation}
By Lemma \ref{lem:H-tri.lift}, 
\begin{equation}\label{equ:red5}
  (\N_A(Q)_{{\varphi'}^{\hH}}, U_{\eta,Q}, \varphi')_{\hH} \geqslant_{\ast} \big(\N_A(F)_{\zeta^{\hH}}, \N_{U_{\eta,Q}}(F), \zeta\big)_{\hH}.
\end{equation}
Applying Lemma \ref{lem:H-tri.tran.} to \eqref{equ:red4} and \eqref{equ:red5},
\begin{equation}\label{equ:red6}
  (A_{\varphi^{\hH}}, G_{\eta}, \varphi)_{\hH} \geqslant_{\ast} \big(\N_A(F)_{\zeta^{\hH}}, \N_{U_{\eta,Q}}(F), \zeta\big)_{\hH}.
\end{equation}

Recall that $\N_{U_{\eta,Q}}(F)$ is the stabilizer of $\delta^{\star_F}$ in $\N_G(F)$, where $\delta^{\star_F}\in\IBr(C_F)$ is the inflation of $(\overline{\delta})^{\star_{\overline{F}}}$ to $C_F$, and $\varrho = \zeta^{\N_G(F)}$.
In fact, $F \in \Rad(G),\varrho\in\IBr(\N_G(F)\which \lambda,|F|)$ and $(F, \varrho)\in f(\psi)$.
Since $Q = F \cap K$ and $K$ is normal, we have 
\begin{equation}\label{equ:red8}
  G_{\eta} \cap \N_A(F) = \N_{G_{\eta}}(F) = \N_{G_{\eta}}(F, Q) = \N_{U_{\eta,Q}}(F).
\end{equation}
We now show that $\big(\hH \times \N_A(F)\big)_{\varrho} = \big(\hH \times \N_A(F)\big)_{\zeta} \N_G(F)$.
Let $a \in \big(\hH \times \N_A(F)\big)_{\varrho}$. Since $C_F \zg \N_A(F)$, both $\delta^{\star_F}$ and $(\delta^{\star_F})^a$ are irreducible constituents of $\varrho_{C_F}$.
By Clifford theory, there exists $x \in \N_G(F)$ such that $(\delta^{\star_F})^{a x} = \delta^{\star_F}$.
Since $\zeta^{a x} \in \IBr(\N_{U_{\eta,Q}}(F) \which \delta^{\star_F})$ is also the Clifford correspondent of $\varrho$, we have $\zeta^{a x} = \zeta$. 
(Note that $\N_{U_{\eta,Q}}(F)^{a x} = \N_{U_{\eta,Q}}(F)$ because $\N_G(F) \zg \N_A(F)$ and $\N_{U_{\eta,Q}}(F)$ is the stabilizer of $\delta^{\star_F}$ in $\N_G(F)$.)
Hence, 
\begin{equation}\label{equ:red7}
  \big(\hH \times \N_A(F)\big)_{\delta} = \big(\hH \times \N_A(F)\big)_{\varrho} \N_G(F).
\end{equation}

Since $f$ is $(\hH \times A)_{\lambda}$-equivariant bijection, and since $ (F,\varrho)= f(\psi)$, it follows as in the proof of Proposition \ref{prop:lieover2} that 
\[
A_{\psi^{\hH}} = G \N_A(F)_{\varrho^{\hH}} \quad \text{and} \quad \big(\hH \times \N_A(F)\big)_{\psi} = \big(\hH \times \N_A(F)\big)_{\varrho}.
\]
By \eqref{equ:red6}, \eqref{equ:red8}, \eqref{equ:red7}, and Theorem \ref{thm:H-tri(ind)},
\[
(A_{\psi^{\hH}}, G, \psi)_{\hH} \geqslant_{\ast} (\N_A(F)_{\varrho^{\hH}}, \N_G(F), \varrho)_{\hH}.
\]
Note that $A_{\varphi^{\hH}} \subseteq A_{\psi^\hH}$ and $\N_A(F)_{\zeta^{\hH}} \subseteq \N_A(F)_{\varrho^{\hH}}$.
This completes the proof of the theorem.
\end{proof}

\section{The inductive BNAW condition for finite non-abelian simple groups
of Lie type  at their defing characteristic}

Let $G$ be a finite group.
Recall that a weight of $G$ can be also defined as a pair $(Q,\varphi)$, where $\varphi\in\Irr(\N_G(Q)/Q)$ with $\varphi(1)_p=|\N_G(Q)/Q|_p$.
In the following, when mentioning a weight, we often mean this version.
Note that $\varphi^{\circ}$, restriction of $\varphi$ to all $p'$-elements of $\N_G(Q)/Q$, is a projective irreducible Brauer character of $\N_G(Q)/Q$. 
In this way, we can identify $(Q,\varphi)$ with $(Q,\varphi^{\circ})$. We also regard $\varphi$ as an irreducible character of $\N_G(Q).$

\begin{prop}\label{def:ind-GAW-cond-2}
  Let $L$ be a finite non-abelian  simple group of order divisible by $p$ and $S$ be the
universal $p'$-covering group of $L$. 
Assume that the following conditions hold.
  \begin{enumerate}
    \item There is an $\hH\times \Aut(S)$-equivariant bijection \[\Omega: \IBr(S)\to\cW(S)\] which preserves blocks. 
    \item For any $\eta\in\IBr(S)$, we let $\overline{(Q,\varphi)}=\Omega(\eta)$. 
        Let $Z$ be an $\Aut(S)_{\eta^{\hH}}$-invariant subgroup of $\Z(S)\cap\ker\eta\cap\ker \varphi^{\circ}$.
      Write $\bar S=S/Z$, $\bar Q=QZ/Z,$ and we also regard $\eta$ as a character of $\bar{S}$ and $\varphi^{\circ}$ as a character of $\N_{\bar S }(\bar Q)$.
    Then $\bar S$ can be normally embedded into a finite group $\widetilde{A},$ and there exists a normal subgroup $A$ of $\widetilde{A}$ containing $\bar S\C_{\widetilde{A}}(\bar S)$, and characters $\tilde\eta\in\IBr(A)$ and $\tilde\varphi\in\IBr(\N_A(\bar Q))$ such that  the following conditions hold.
        \begin{enumerate}
          \item $\widetilde{A}/\C_{\widetilde{A}}(\bar S)\cong \Aut(\bar S)_{\eta^{\hH}}$ and $A/\C_{\widetilde{A}}(\bar S)\cong \Aut(\bar S)_{\eta}.$ Moreover, $\C_{\widetilde{A}}(\bar S)=\Z(A)$. 
          \item $\tilde\eta$ is an extension of $\eta$ and $\tilde\varphi$ is an extension of $\varphi^{\circ}$ such that $\IBr(\Z(A)\which \tilde\eta)=\IBr(\Z(A)\which  \tilde\varphi)$.
          \item For any $J\in\mS(A,\bar S),$ we have $\bl(\tilde{\varphi}_{\N_J(\bar Q)})^J=\bl(\tilde\eta_J).$
          \item For any $a\in(\hH\times \N_{\widetilde{A}}(\bar Q))_\varphi$, let $\tilde\eta^{a}=\mu_a\tilde\eta$ and $\tilde\varphi^a=\mu'_a\tilde\varphi$, where $\mu_a,\mu'_a$ are linear Brauer characters of $A/\bar S=\N_A(\bar Q)/\N_{\bar S}(\bar Q)$, then we have $\mu_a=\mu'_a.$
        \end{enumerate}
  \end{enumerate}
Then the inductive BNAW condition holds for $L$ at the  prime $p$.
\end{prop}

\begin{proof}
We only need to make some adjustment to the proof of \cite[Proposition 9.1]{FFZ}. 
Fix any $\eta\in\IBr(S)$ and $\overline{(Q,\varphi)}=\Omega(\eta)$.
Let $\widetilde{G}=S\semi \Aut(S)_{\eta^{\hH}}.$
We are left to prove that
\begin{equation}\label{equ:p9.1}
(\widetilde{G},S,\eta)_{\hH}\geqslant_b (\N_{\widetilde{G}}(Q),\N_S(Q),\varphi^{\circ})_{\hH}.
\end{equation}
By (1) we have $\widetilde{G}=S\N_{\widetilde{G}}(Q)$ and $(\hH\times\N_{\widetilde{G}}(Q))_{\eta}= (\hH\times\N_{\widetilde{G}}(Q))_{\varphi^{\circ}}.$ 
Note that $\overline{\N_S(Q)}=\N_{\bar S}(\bar Q)$ (see the proof of  \cite[Proposition 9.1]{FFZ}).
A direct computation shows that \cite[Propsition 2.3]{NS14} holds for Brauer characters.
Thus condition (2) tells that 
\begin{equation*}
  (\widetilde{A},\bar{S},\eta)_{\hH} \geqslant_b(\N_{\widetilde{A}}(\bar{Q}), \N_{\bar{S}}(\bar{Q}),\varphi^{\circ})_{\hH}.
\end{equation*}
Since $\N_A(\bar Q)$ and $\N_{\widetilde{G}/Z}(\bar Q)$ afford the same automorphism group of $\bar S$ (see the proof of  \cite[Proposition 9.1]{FFZ} for details), by Theorem \ref{thm:butterfly}, we have $$(\widetilde{G}/Z,\bar{S},\eta)_{\hH} \geqslant_b(\N_{\widetilde{G}}({Q})/Z,\N_{\bar{S}}(\bar Q),\varphi^{\circ})_{\hH}.$$  
Then \ref{equ:p9.1} follows from Lemma \ref{lem:H-tri.lift}. 
This completes proof of the proposition.
\end{proof}  

By \cite[Corollary 10.4]{FFZ}, the inductive NAW condition holds for simple groups with a cyclic Sylow 
$p$-subgroup and a cyclic outer automorphism group. We show that, in this situation, the inductive BNAW condition also holds. We need to strengthen \cite[Proposition 10.3]{FFZ} slightly.
Let $G$ be a finite group and $Q$ be $p$-subgroup of $G$. We denote by $\IBrd(G \which Q)$ the set of irreducible
Brauer characters of $G$ lying in blocks with defect group $Q$.

\begin{coro}\label{coro:10.4}
Let $L$ be a finite non-abelian  simple group of order divisible by $p$ and  $S$ be the universal $p'$-covering group of $L$. Suppose that there exists a finite group $G$ such that $S\zg G$ and $G$ induces all automorphisms on $S$ by conjugation.
Assume that $G/S$ is cyclic, $\C_G(S)=\Z(G),$ and $S$ has a cyclic Sylow $p$-subgroup $P$. Then the inductive BNAW condition holds for $L$ at the prime $p$.
\end{coro}
\begin{proof}
  Keep the notation in the proof of \cite[Corollary 10.4]{FFZ}.
  We have proved that $(\tilde{\eta},\tilde\varphi)$ gives \[(\widetilde A,S,\eta)_{\hH}\geqslant_{c} (\N_{\widetilde A}(Q),\N_S(Q),\varphi)_{\hH}.\]
  Since the centralizer in $\widetilde{A}_{\eta}$ of a defect group of $\bl(\varphi)$ is contained in $\N_{\widetilde A}(Q)_{\varphi}.$
  Thus by \cite[Propsition 2.3]{NS14}, we only need to prove that 
  \begin{equation}\label{equ:10.1}
    \bl(\tilde\varphi_{\N_J(Q)})^J=\bl(\tilde{\eta}_J)
  \end{equation} for any $J\in\mS(A,S),$ where $A=\widetilde{A}_{\eta}.$
  Let $J_1=J\cap A_1.$
  Since the characters
  $\tilde{\eta}_{J_1}$ and $\tilde\varphi_{\N_{J_1}(Q)}$ correspond to each other through the bijection constructed in \cite[Proposition 10.3]{FFZ},   we have
  \begin{equation}\label{equ:10.2}
    \bl(\tilde\varphi_{\N_{J_1}(Q)})^{J_1}=\bl(\tilde{\eta}_{J_1}).
  \end{equation}
Then \ref{equ:10.1} follows from \ref{equ:10.2} and \cite[Lemma 2.3]{KS15}, as $\bl(\tilde\eta_J)$ is the unique block of $J$ covering the block $\bl(\tilde\eta_{J_1})$ of $J_1.$
\end{proof}

  \begin{proof}[Proof of Theorem D]
    Let $L$ be a finite non-abelian simple group of Lie type  in characteristic $p$, and $S$ be the universal $p'$-covering group of $L$.
    
      First, we show that the simple group $L\in\{\textup{Sp}_4(2)'\cong A_6, \textup G_2(2)'\cong\textup{SU}_3(3), {}^2\textup F_4(2)'\}$ satisfies the inductive BNAW condition at the prime 2.
    Recall that $L$ has trivial Schur multiplier (so we set $S:=L$), and $\textup{Out}(S)$ is of order 2 or isomorphic to the Klein four group.
    In the proof of \cite[Prop.~9.3]{FFZ}, a $\mathcal H\times\Aut(S)$-equivariant bijection $\Omega: \IBr(S)\to\cW(S)$ is already established, and it suffices to show that it preserves blocks.
    By the construction of $\Omega$, we see that it preserves the blocks of defect zero. 
    On the other hand, in  the proof of \cite[Prop.~9.3]{FFZ} it is shown that the actions of $\mathcal H$ on both irreducible Brauer characters and (conjugacy classes of) weights in blocks of positive defect are trivial.
    Hence we only need to prove that there exists an $\Aut(S)$-equivariant bijection $\Omega: \IBr(S)\to\cW(S)$. This has been obtained in the proof of inductive BAW condition in \cite[p.215]{Spa13}. So the inductive BNAW condition holds for the simple groups $\textup{Sp}_4(2)'$, $\textup G_2(2)'$ and ${}^2\textup F_4(2)'$ at the prime 2.
    
    The group $\textup{SL}_2(8)$ has cyclic Sylow 3-subgroups, and thus the inductive BNAW condition holds for the group $\textup{SL}_2(8)$ at the prime 3 by Corollary~\ref{coro:10.4}.
    Therefore, we may assume that $L\notin\{\textup{Sp}_4(2)', \textup G_2(2)', {}^2\textup F_4(2)'\}$ when $p=2$, and $L\ncong {}^2G_2(3)'\cong\textup{SL}_2(8)$ when $p=3$.
    For every other situation the exceptional part of the Schur multiplier of $L$ is a $p$-group (see, e.g., \cite[Table 6.1.3]{GLS98}).
    Hence $S=\mathcal G^F$, for some simply-connected simple algebraic group $\mathcal G$ defined over $\overline{\mathbb F}_p$ and some Steinberg map $F:\mathcal G\to\mathcal G$.
    
    In the proof of \cite[Thm.~C]{FFZ}, we have proved that the inductive NAW condition holds for $L$, and thus by construction, it suffices to show that the condition (2.c) of Proposition~\ref{def:ind-GAW-cond-2} holds.
    Now we slightly modify the bijection between irreducible Brauer characters and conjugacy classes of weights in the proof of \cite[Thm.~C]{FFZ}, swapping the principal character and the Steinberg character (following \cite[p.~216]{Spa13}).
    Note that the correspondences between irreducible Brauer characters and conjugacy classes of weights, and the extensions of Brauer characters and weight characters constructed in the proofs of \cite[Thm.~C]{FFZ} and \cite[Thm.~C]{Spa13} coincide.
    Therefore, by the proof of \cite[Thm.~C]{Spa13}, the block inductions in the condition (2.c) of Proposition~\ref{def:ind-GAW-cond-2} holds, which completes the proof.
    \end{proof}

\bigskip
\noindent\textbf{Acknowledgements}
\bigskip

The authors would like to express their sincere gratitude to Xueqin Hu for his invaluable suggestions on the proof of Theorem \ref{thm:DGN1}.
The authors are supported by the NSFC [Grant Nos. 12350710787, 12571016].


\begin{thebibliography}{99}

\bibitem{Al87}
J.\,L.~Alperin, Weights for finite groups. In: \emph{The Arcata Conference on Representations of Finite Groups, Arcata, Calif. (1986), Part I}. Proc. Sympos. Pure Math., vol. \textbf{47}, pp. 369--379, Am. Math. Soc., Providence, 1987.

\bibitem{Dade}
E.\,C.~Dade, A correspondence of characters. In: {\itshape The Santa Cruz Conference on Finite Groups}, Univ. California, Santa
Cruz, CA, 1979, in: Proc. Sympos. Pure Math., vol. {\bf 37}, Amer. Math. Soc., Providence, RI, 1980, pp. 401–403.

\bibitem{Duhu25}
Y. Du, X. Huang, 
The blockwise Galois Alperin weight conjecture for symmetric and alternating groups.
\emph{Math. Z. \bf 311} (2025), no. 1, Paper No. 3, 19 pp.

\bibitem{Fe25}
Z. Feng, Character triples and weights. \emph{J. Algebra \bf 684} (2025), 690--725.

\bibitem{FFZ}
Z. Feng, Q. Fu, Y. Zhou, A reduction theorem for the Navarro Alperin wight conjecture. arXiv:2312.02594.

\bibitem{FZ22}
Z. Feng, J. Zhang, Alperin weight conjecture and related developments.
  \emph{Bull. Math. Sci. \bf 12} (2022), 2230005.

\bibitem{Fu24}
Q. Fu, An equivariant bijection of irreducible Brauer characters above the Dade--Glauberman--Nagao correspondence. \emph{J. Algebra \bf 689} (2026), 710--738.

\bibitem{GLS98}
D. Gorenstein, R. Lyons, R. Solomon, \emph{The Classification of the Finite Simple Groups, Number 3}. Math. Surveys Monogr., {vol. \bf 40}, American Mathematical Society, Providence, RI, 1998.

\bibitem{IMN07}
I.\,M. Issacs, G. Malle, G. Navarro, A reduction theorem for the McKay conjecture. {\itshape Invent. Math.} \textbf{170} (2007), 33--101.

\bibitem{KS15}
S. Koshitani, B. Sp\"ath, Clifford theory of characters in induced blocks. {\itshape Proc. Amer. Math. Soc.} \textbf{143} (2015), 3687--3702.

\bibitem{Linckelmann}
M. Linckelmann, {\itshape The Block Theory of Finite Group Algebras, volume II}.  London Math. Soc. Stud. Texts, {vol. \bf~92}, Cambridge University Press, Cambridge, 2018.

\bibitem{MRR23}
J.\,M.~Martínez, N.~Rizo, D.~Rossi, The Alperin weight conjecture and the Glauberman correspondence via character triples.
{\itshape Algebra Number Theory} \textbf{20} (2026), no. 2, 333–382. 

\bibitem{Mc71}
J. McKay, A new invariant for simple groups. {\itshape Notices. Amer. Math. Soc.} \textbf{18} (1971), 397. 

\bibitem{Mc72}
J. McKay, Irreducible representations of odd degree. {\itshape J. Algebra} \textbf{20} (1972), 416--418.

\bibitem{NagTsu}
H. Nagao, Y. Tsushima, {\itshape Representations of Finite Groups}. Academic Press, Boston, 1989.

\bibitem{Nav:char.}
G. Navarro, {\itshape Characters and Blocks of Finite Groups}. London Mathematical Society Lecture Note Series, {vol. \bf~250}. Cambridge University Press, Cambridge, 1998.

\bibitem{Na04}
G. Navarro, The Mckay conjecture and Galois automorphisms. {\itshape Ann. of Math. (2)} \textbf{160} (2004), 1129--1140.

\bibitem{Navarro:McKay}
G. Navarro, {\itshape Character Theory and the McKay Conjecture}. Cambridge Studies in Advanced Mathematics, vol. {\bf 175}, Cambridge University Press, Cambridge, 2018. 

\bibitem{NS14}
G. Navarro, B. Sp\"ath, On Brauer's height zero conjecture. {\itshape J. Eur. Math. Soc.} \textbf{16} (2014), 695--747.

\bibitem{NSV20}
G. Navarro, B. Sp\"ath, C. Vallejo, A reduction theorem for the Galois--Mckay conjecture. {\itshape Trans. Amer. Math. Soc.} \textbf{373} (2020), 6157--6183.

\bibitem{NT11}
G. Navarro, P.\,H. Tiep, A reduction theorem for the Alperin weight conjecture. {\itshape Invent. Math.} \textbf{184} (2011), 529--565.



\bibitem{Ro23}
D. Rossi, The McKay conjecture and central isomorphic character triples. {\itshape J. Algebra} \textbf{618} (2023), 42--55.

\bibitem{Ro25}
D. Rossi, A reduction theorem for the Character Triple Conjecture.
\emph{Proc. Lond. Math. Soc.~(3)~\bf 131} (2025), No. e70090.

\bibitem{Spa13}
B. Sp\"ath, A reduction theorem for the blockwise Alperin weight conjecture. {\it J. Group Theory} \textbf{16} (2013), 159--220.

\bibitem{Spa17}
B. Sp\"ath, A reduction theorem for Dade's projective conjecture. {\itshape J. Eur. Math. Soc.}  \textbf{19} (2017), 1071--1126.

\bibitem{Spath}
B. Sp\"ath, Inductive conditions for counting conjectures via character triples. In: \emph{Representation Theory -- Current Trends and Perspectives}. EMS Ser. Congr. Rep., Eur. Math. Soc., Z\"urich, 2017, pp.~665--680.

\bibitem{SV16} 
B. Sp\"ath, C. Vallejo, Brauer characters and coprime action. {\itshape J. Algebra} \textbf{457} (2016), 276--311.

\bibitem{Thevenaz}
J. Thévenaz, 
{\itshape $G$-algebras and Modular Representation Theory}. The Clarendon Press, Oxford University Press, New York,  1995.

\bibitem{Tu94}
A. Turull, Clifford theory with Schur indices. {\itshape J. Algebra} \textbf{170} (1994), 661--677.

\bibitem{Tu08}
A. Turull, Above the Glauberman correspondence. {\itshape Adv. Math.} \textbf{217} (2008), 2170--2205.

\bibitem{Tu09:BC}
A. Turull, The Brauer--Clifford group. {\itshape J. Algebra} \textbf{321} (2009), 3620--3642.

\bibitem{Tu09:FM}
A. Turull, Brauer--Clifford equivalence of full matrix algebras. {\itshape J. Algebra} \textbf{321} (2009), 3643--3658.

\bibitem{Tu14}
A. Turull, The strengthened Alperin weight conjecture for $p$-solvable groups. {\itshape J. Algebra} \textbf{398} (2014), 469--480.

\end{thebibliography}
\end{document}